\documentclass{amsart}

\usepackage{amsmath,amssymb}
\usepackage{amsthm}
\usepackage{comment}
\usepackage{mathtools}
\usepackage{graphics}
\usepackage{wrapfig}
\usepackage{hyperref}
\usepackage{subcaption}
\usepackage{cite}
\usepackage{tikz}

\usepackage{algorithm}
\usepackage{algpseudocode}
\usepackage{booktabs}
\floatname{algorithm}{Algorithm}

\newtheorem{theorem}{Theorem}[section]
\newtheorem{lemma}[theorem]{Lemma}
\newtheorem{proposition}[theorem]{Proposition}
\newtheorem{cor}[theorem]{Corollary}
\newtheorem{definition}[theorem]{Definition}
\newtheorem{ex}[theorem]{Example}
\newtheorem{remark}[theorem]{Remark}

\DeclareMathOperator{\Def}{Def}
\DeclareMathOperator{\Diam}{Diam}
\graphicspath{ {./images/} }

\numberwithin{equation}{section}

\newcommand{\Id}{\mathrm{Id}}

\newcommand{\bF}{\mathbb{F}}

\newcommand{\bP}{\mathbb{P}}

\newcommand{\cP}{\mathcal{P}}
\newcommand{\Z}{\mathbb{Z}}

\newcommand{\R}{\mathbb{R}}

\newcommand{\cU}{\mathcal{U}}

\newcommand{\cM}{\mathcal{M}}
\newcommand{\cS}{\mathcal{S}}
\newcommand{\cF}{\mathcal{F}}
\newcommand{\cI}{\mathcal{I}}
\newcommand{\cA}{\mathcal{A}}

\newcommand{\gr}{\mathrm{gr}}

\newcommand{\cL}{\mathcal{L}}
\newcommand{\cZ}{\mathcal{Z}}

\begin{document}
\title{A logifold structure for measure space}
\author{Inkee Jung and Siu-Cheong Lau}

\begin{abstract}
In this paper, we develop a geometric formulation of datasets. The key novel idea is to formulate a dataset to be a fuzzy topological measure space as a global object, and equip the space with an atlas of local charts using graphs of fuzzy linear logical functions. We call such a space to be a logifold. In applications, the charts are constructed by machine learning with neural network models. We implement the logifold formulation to find fuzzy domains of a dataset and to improve accuracy in data classification problems.
\end{abstract}

\maketitle{}

\section{Introduction}
In geometry and topology,  the manifold approach dated back to Riemann uses open subsets in $\R^n$ as local models to build a space.  Such a local-to-global principle is central to geometry and has achieved extremely exciting breakthroughs in modeling spacetime by Einstein's theory of relativity.  

In recent years, the rapid development of data science brings immense interest to datasets that are `wilder' than typical spaces that are well-studied in geometry and topology.  
Taming the wild is a central theme in the development of mathematics.  Advances in computational tools have helped to expand the realm of mathematics in the history.
For instance, it took many years in human history to recognize the irrational number $\pi$ and approximate it by rational numbers.  In this regard, we consider machine learning by neural network models as a modern tool to find expressions of a `wild space' (for instance a dataset in real life) as a union (limit) of fuzzy spaces expressed by finite formulae.

Let $X$ be a topological space, $B_X$ the corresponding Borel $\sigma$-algebra, and $\mu$ a measure on $(X,B_X)$.  We understand a dataset as a \emph{fuzzy topological measure space} in nature. To work with such a complicated space, we would like to have local charts that admit finite mathematical expressions and has logical interpretations.  This plays the role of local coordinate systems for a measure space. In our definition, a local chart $U \subset X$ is required to be have positive measure.
Moreover, to avoid triviality and requiring too many charts to cover the whole space, we may further fix $\epsilon>0$ and require that $\mu(X-U)<\epsilon$.
Such a condition disallows $U$ to be too simple, such as a tiny ball around a point in a dataset.
This resembles the Zariski-open condition in algebraic geometry.

In place of open subsets of $\R^n$, we formulate `local charts' that are closely related to neural network models and have logic gate interpretations.  Neural network models provide a surprisingly successful tool to find mathematical expressions that approximate a dataset.  Non-differentiable or even discontinuous functions analogous to logic gate operations are frequently used in network models.  It provides an important class of non-smooth and even discontinuous functions to study a space.
 
We take classification problems as the main motivation in this paper.
For this, we consider the graph of a function $f: D \to T$, where $D$ is a measurable subset in $\R^n$ (with the standard Lebesgue measure) and $T$ is a finite set (with the discrete topology).
The graph $\gr(f) \subset D \times T$ is equipped with the push-forward measure by $D \to \gr(f)$.  

We will use the graphs of linear logical functions $f: D \to T$ explained below as local models.
A chart is of the form $(U,\Phi)$, where $U \subset X$ is a measurable subset which satisfies $\mu(X-U)<\epsilon$, and $\Phi: U \to \gr(f)$ is a measure-preserving homeomorphism.
We define a \emph{linear logifold} to be a pair $(X, \cU)$ where $\cU$ is a collection of charts $(U_i,\Phi_i)$ such that $\mu(X-\bigcup_i U_i) = 0$.

The definition of linear logical functions is motivated from neural network models and has a logic gate interpretation.  A network model consists of a directed graph, whose arrows are equipped with linear functions and vertices are equipped with non-linear functions, which are typically ReLu or sigmoid functions in middle layers, and are sigmoid or softmax functions in the last layer.  Note that sigmoid and softmax functions are smoothings of the discrete-valued step function and the index-max function respectively.  Such smoothings are useful to describe fuzziness of data.

From this perspective, step and index-max functions are the non-fuzzy (or called classical) limit of sigmoid and softmax functions respectively.  We will take such a limit first, and come back to fuzziness in a later stage. This means we replace all sigmoid and softmax functions in a neural network model by step and index-max functions. We will show that such a neural network is equivalent to a linear logical graph:
at each node of the directed graph, there is a system of $N$ linear inequalities (on the input Euclidean domain $D \subset \R^n$) that produces $2^N$ possible Boolean outcomes for an input element, which determine the next node that this element will be passed to.  
We call the resulting function to be a linear logical function.

We prove that linear logical functions can approximate any given measurable function $f:D \to T$, where $D$ is a measurable subset of $\R^n$ with $\mu(D)<\infty$.  This provides a theoretical basis of using these functions in modeling.

\begin{theorem}[Universal approximation theorem by linear logical functions]
	Let $f:D\to T$ be a measurable function whose domain $D \subset \R^n$ is of finite Lebesgue measure, and suppose that its target set $T$ is finite.  For any $\epsilon>0$, there exists a linear logical function $L$ and a measurable set $E \subset \R^n$ of the Lebesgue measure less than $\epsilon$ such that $L|_{D-E} \equiv f|_{D-E}$.
\end{theorem}
By taking the limit $\epsilon \to 0$, the above theorem finds a linear logifold structure on the graph of a measurable function.  In reality, $\epsilon$ reflects the error of a network in modeling a dataset.  

It turns out that for $D=\R^n$, linear logical functions $\R^n \to T$ (where $T$ is identified with a finite subset of $\R$) are equivalent to semilinear functions $\R^n \to T$, whose graphs are semilinear sets defined by linear equations and inequalities \cite{vandenDries}.
Semilinear sets provide the simplest class of definable sets of so-called o-minimal structures \cite{vandenDries}, which are closely related to model theory in mathematical logic. 
O-minimal structures made an axiomatic development of Grothendieck's idea of finding tame spaces that exclude wild topology.  On one hand, definable sets have a finite expression which is crucial to the predictive power and interpretability in applications.
On the other hand, our setup using measurable sets $D \subset \R^n$ provides a larger flexibility for modeling data.

Compared to more traditional approximation methods such as Fourier series, there are reasons why linear logical functions are preferred in many situations for data.  When the problem is discrete in nature (for instance the target set $T$ is finite), it is simple and natural to take the most basic kinds of discrete-valued functions as building blocks, namely step functions formed by linear inequalities.  These basic functions are composed to form networks which are supported by current computational technology.  Moreover, such discrete-valued functions have fuzzy and quantum deformations which have rich meanings in mathematics and physics. 

\begin{figure}[ht]
	\includegraphics[scale=0.6]{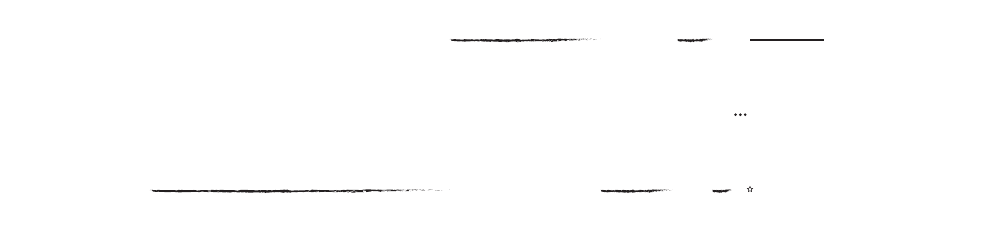}
	\caption{An example of a logifold.  The graph jumps over values $0$ and $1$ infinitely in left-approaching to the point marked by a star (and the length of each interval is halved).  This is covered by infinitely many charts of  linear logical functions, each of which has only finitely many jumps. Moreover, the base is a measurable subset of $\R$ (which is hard to depict and not shown in the picture).}
	\label{fig:limit-logifold}
\end{figure}

Now let us address \emph{fuzziness}, another important feature of a dataset besides discontinuities.  In practice, there is always an ambiguity in determining whether a point belongs to a dataset. This is described as a fuzzy space $(X,\cP)$, where $X$ a topological measure space and $\cP: X \to (0,1]$ is a continuous measurable function that encodes the probability of whether a given point of $X$ belongs to the fuzzy space under consideration. Here, we require $\cP>0$ on purpose. While points of zero probability of belonging may be adjoined to $X$ to that the union gets simplified (for instance $X$ may be embedded into $\R^n$ where points in $\R^n - X$ has zero probability of belonging), they are auxiliary and have no intrinsic meaning. 

To be useful, we need a finite mathematical expression (or approximation) for $\cP$. This is where neural network models enter into the description. A neural network model for a classification problem that has the softmax function in the last layer gives a function $f = (f_1,\ldots,f_d):\R^n \to S$ where $S$ is the standard simplex $\{\sum_{i=1}^d y_i = 1\} \subset \R^d$. This gives a fuzzy space 
$$(\R^n \times T,\cP)$$ 
where $T=\{1,\ldots,d\}$ and $\cP(p,t) := f_t(p)$. As we have explained above, in the non-fuzzy limit sigmoid and softmax functions are replaced by their classical counterparts of step and index-max functions respectively, and we obtain $f^{\textrm{classical}}: \R^n \to T$ and the subset $\{(p,t):f^{\textrm{classical}}(p)=t\} \subset \R^n \times T$ as the classical limit. 

However, the ambient space $\R^n$ is not intrinsic: for instance, in the context of images,  the dimension $n$ gets bigger if we take images with higher resolutions, even though the objects under concern remain the same. Thus, like in manifold theory the target space is taken to be a topological space rather than $\R^n$, our theory take a topological measure space $X$ in place of $\R^n$ (or $\R^n \times T$). $\R^n \times T$ (for various possible $n$) is taken as an auxiliary ambient space that contains (fuzzy) measurable subsets that serve as charts to describe a dataset $(X,\cP)$. 

Generally, we formulate fuzzy linear logical functions (Definition \ref{def:fuzlinlogfcn}) and fuzzy linear logifolds (Definition \ref{def:fuzlogifold}). 
A fuzzy logical graph is 
a directed graph $G$  whose each vertex of $G$ is equipped with a state space and each arrow is equipped with a continuous map between the state spaces.  The walk on the graph (determined by inequalities) depends on the fuzzy propagation of the internal state spaces.

\begin{figure}[ht]
	\includegraphics[scale=0.6]{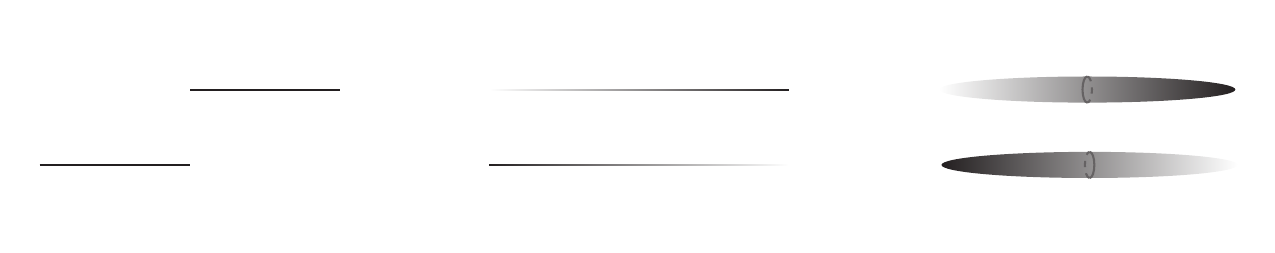}
	\caption{The left hand side shows a simple example of a logifold.  It is the graph of the step function $[-1,1] \to \{0,1\}$.  The figure in the middle  shows a fuzzy deformation of it, which is a fuzzy subset in $[-1,1]\times\{0,1\}$.  The right hand side shows the graph of probability distribution of a quantum observation, which consists of the maps $\frac{|z_0|^2}{|z_0|^2 + |z_1|^2}$ and $\frac{|z_1|^2}{|z_0|^2 + |z_1|^2}$ from the state space $\bP^1$ to $[0,1]$.}
	\label{fig:discrete-fuzzy-quantum}
\end{figure}

For readers who are more computationally oriented, they can first directly go to read Section \ref{sec:algo} where we describe the implementation of the logifold theory in algorithms. 
Our logifold formulation of a dataset can be understood as a geometric theory for ensemble learning and the method of Mixture of Experts. 

Ensemble learning utilizes multiple trained models to make a decision or prediction, see for instance \cite{EnsembleML}.
Ensemble machine learning achieves improvement in classification problems, see for instance \cite{KNNensembleCIFAR} and \cite{EfficientNet}. In the method of Mixture of Experts \cite{AdaptiveMixtureofLocalExperts}, several expert models $E_j$ are employed, and there is also a gating function $G_j(x)$. The final outcome is given by the total $\sum_j G_j(x) E_j(x)$. This idea of using `experts' to describe dataset is similar to the formulation of a fuzzy logifold. On the other hand, motivated from manifold theory, we formulate \emph{universal mathematical structures that are common to datasets, namely the global intrinsic structure of a fuzzy topological measure space, and local logical structures among data points expresssed by graphs of fuzzy logical functions}.

We also propose refinement and enhancement in algorithms.
The key new ingredient in our implementation is the fuzzy domain of each model.
A trained model typically does not have a perfect accuracy rate and performs well only on a subset of data points, or for a subset of target classes. A major step here is to find and record the domain of each model where it works well.

We restrict the domain of a model to $D \times T'$ for some subsets $D\subset \R^n$ and $T'\subset T$. 
Certainty scores are used to determine the (fuzzy) scopes of models, with the softmax function in the last layer of each node providing these scores.
For each network model, we restrict the input domain to the subset of data that has certainty scores higher than a threshold.
The final threshold for testing data can be set based on expected accuracy obtained from fuzzy domains measured on validation data. We successfully improve the accuracy (by 5\% to more than 20\% depending on the types of experiments) using fuzzy domains and our refined voting method \cite{Jung-Lau}.

The structure of this paper is as follows. First, we formulate fuzzy linear logical functions in Section \ref{sec:functions}. Next, we establish relations with semilinear functions in Section \ref{sec:semilinear}, prove the universal approximation theorem for linear logical functions in Section \ref{sec:universal}, and define fuzzy linear logifolds in Section \ref{sec:logifold}. We provide a detailed description of algorithmic implementation of logifolds in Section \ref{sec:algo}.

\section{Linear logical functions and their fuzzy analogs} \label{sec:functions}

Given a subset of $\R^n$, one would like to describe it as the zero locus, the image, or the graph of a function in a certain type.  
In analysis, we typically think of continuous/smooth/analytic functions.  However, when the domain is not open, smoothness may not be the most relevant condition. 

The success of network models has taught us a new kind of functions that are surprisingly powerful in describing datasets.  Here, we formulate them using directed graphs and call them \emph{linear logical functions}.
The functions offer three distinctive advantages.
First, they have the advantage of being logically interpretable in theory.  Second, they are close analogues of quantum processes.  Namely, they are made up of linear functions and certain non-linear activation functions, which are analogous to unitary evolution and quantum measurements.
Finally, it is natural to add fuzziness to these functions and hence they are better adapted to describe statistical data.

\subsection{Linear logical functions and their graphs}

We consider functions $D \to T$ for $D \subset \R^n$ and a finite set $T$ constructed from a graph as follows. 
Let $G$ be a finite directed graph that has no oriented cycle and has exactly one source vertex which has no incoming arrow and $|T|$ target vertices.
Each vertex that has more than one outgoing arrows is equipped with an affine linear function $l=(l_1,\ldots,l_k)$ on $\R^n$, where the outgoing arrows at this vertex are one-to-one corresponding to the chambers in $\R^n$ subdivided by the hyperplanes $\{l_i=0\}$.  (For theoretical purpose, we define these chambers to contain some of their boundary strata in a way such that they are disjoint and their union equals $\R^n$.)

\begin{figure}[ht]
	\includegraphics[scale=0.5]{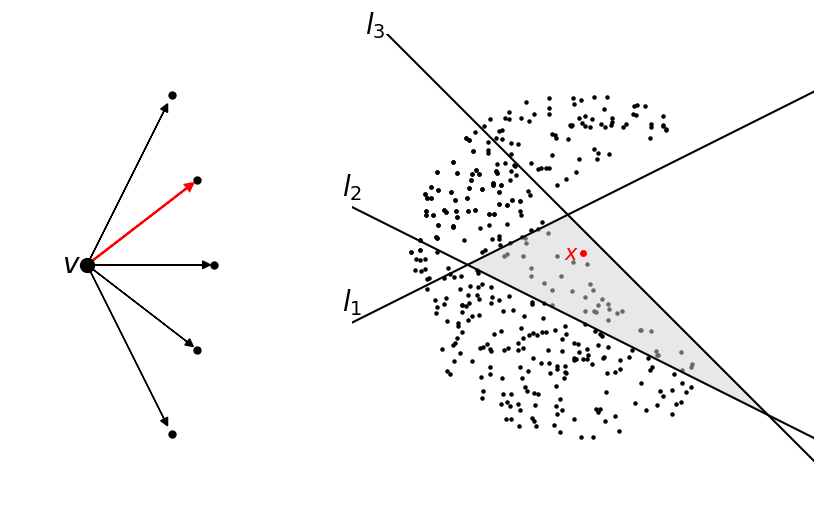}
	\caption{The left side shows a partial directed graph at vertex $v$, with five outgoing arrows. On the right, chambers are formed in $\R^2$ by the affine maps $L_v = (l_1, l_2, l_3)$ defined on $\R^2$. A point $x$ is marked in the chamber defined by  $\{l_1 \leq 0, l_2 \geq 0, l_3 \leq 0\}$. One of the arrows corresponding to the shaded chamber containing $x$ is highlighted in the left diagram.}
	\label{fig:ChamberAndArrow}
\end{figure}

\begin{definition} \label{def:logfcn}
	A linear logical function $f_{G,L}: D \to T$ is a function made in the following way from $(G,L)$,
	where $G$ is a finite directed graph that has no oriented cycle and has exactly one source vertex and $|T|$ target vertices, \[L= \{l_v: v \textrm{ is a vertex with more than one outgoing arrows}\},\]
	$l_v$ are affine linear functions whose chambers in $D$ are one-to-one corresponding to the outgoing arrows of $v$.  $(G,L)$ is called a linear logical graph.
	
	Given $x\in D$, we get a path from the source vertex to one of the target vertices in $G$ as follows.  We start with the source vertex.  At a vertex, if there is only one outgoing arrow, we simply follow that arrow to reach the next vertex.  If there are more than one outgoing arrow, we consider the chambers made by the affine linear function $l$ at that vertex, and pick the outgoing arrow that corresponds to the chamber that $x$ lies in. See Figure \ref{fig:ChamberAndArrow}. Since the graph is finite and has no oriented cycle, we will stop at a target vertex, which is associated to an element $t\in T$.  This defines the function $f_{G,L}$ by setting $f_{G,L}(x)=t$.
\end{definition}

\begin{proposition}
	Consider a feed-forward network model whose activation function at each hidden layer is the step function, and that at the last layer is the index-max function.  The function is of the form
	\[\sigma \circ L_{N} \circ s_{N-1} \circ L_{N-1} \circ \ldots \circ s_1 \circ L_1\]
	where $L_i:\R^{n_{i-1}}\to\R^{n_i}$ are affine linear functions with $n_0=n$, $s_i$ are the entrywise step functions and $\sigma$ is the index-max function.  We make the generic assumption that the hyperplanes defined by $L_i$ for $i=2,\ldots,N$  do not contain $s_{i-1}\circ \ldots \circ L_1(D)$.  Then this is a linear logical function with target $T=\{1,\ldots,n_N\}$ (on any $D\subset \R^n$ where $\R^n$ is the domain of $L_1$).
\end{proposition}	
	
\begin{proof}
	The linear logical graph $(G,L)$ is constructed as follows.  The source vertex $v_0$ is equipped with the affine linear function $L_1$.  Then we make $N$ number of outgoing arrows of $v_0$ (and corresponding vertices) where $N$ is the number of chambers of $L_1$, which are one-to-one corresponding to the possible outcomes of $s_1$ (which form a finite subset of ${\{0,1\}}^{n_1}$).  Then we consider $s_2 \circ L_2$ restricted to this finite set, which also has a finite number of possible outcomes.  This produces exactly one outgoing arrow for each of the vertices in the first layer.  We proceed inductively.  The last layer $\sigma \circ L_{N}$ is similar and has $n_N$ possible outcomes.  Thus we obtain $(G,L)$ as claimed, where $L$ consists of only one affine linear function $L_1$ over the source vertex.  
\end{proof}	

\begin{figure}[ht]
	\includegraphics[scale=0.6]{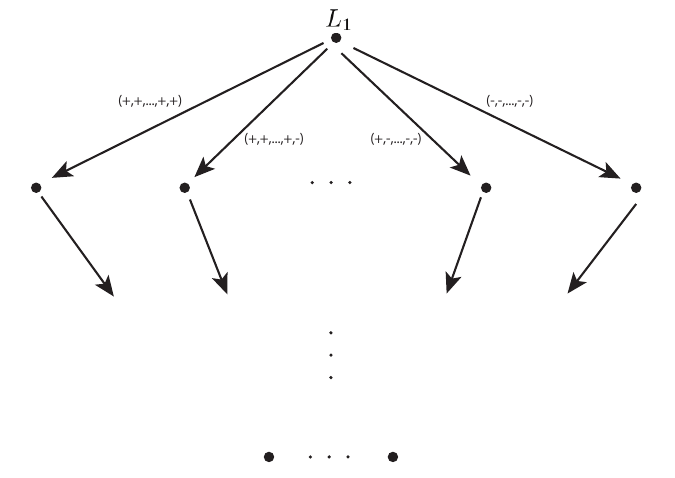}
	\caption{}
	\label{fig:linLogGrStep}
\end{figure}

Figure \ref{fig:linLogGrStep} depicts the logical graph in the above proposition.

\begin{proposition} \label{ex:ReLu}
	Consider a feed-forward network model whose activation function at each hidden layer is the ReLu function, and that at the last layer is the index-max function.  The function takes the form
	  \[\sigma \circ L_{N} \circ r_{N-1} \circ L_{N-1} \circ \ldots \circ r_1 \circ L_1\]  
	where $L_i:\R^{n_{i-1}}\to\R^{n_i}$ are affine linear functions with $n_0=n$, $r_i$ are the entrywise ReLu functions and $\sigma$ is the index-max function.  This is a linear logical function.
\end{proposition}

\begin{proof}
	We construct a linear logical graph $(G,L)$ which produces this function.
The first step is similar to the proof of the above proposition.
Namely, the source vertex $v_0$ is equipped with the affine linear function $L_1$.
Next, we make $N$ number of outgoing arrows of $v_0$ (and corresponding vertices), where $N$ is the number of chambers of $L_1$, which are one-to-one corresponding to the possible outcomes of the sign vector of $r_1$ (which form a finite subset of ${\{0,+\}}^{n_1}$).
Now we consider the next linear function $L_2$.
For each of these vertices in the first layer,  we consider $L_2 \circ r_1 \circ L_1$ restricted to the corresponding chamber, which is a linear function on the original domain $R^n$, and we equip this function to the vertex.
Again, we make a number of outgoing arrows that correspond to the chambers in $R^n$ made by this linear function.
We proceed inductively, and get to the layer of vertices that correspond to the chambers of $L_{N-1} \circ r_{N-2} \circ L_{N-3} \circ \ldots \circ r_1\circ L_1$.
Write $L_N = (l_1,\ldots,l_{n_N})$, and consider $\tilde{L}_N = (l_i-l_j: i\neq j)$.
At each of these vertices, 
\[\tilde{L}_N \circ r_{N-1} \circ L_{N-1} \circ \ldots \circ r_1 \circ L_1\] restricted on the corresponding chamber is a linear function on the original domain $\R^n$, and we equip this function to the vertex and make outgoing arrows corresponding to the chambers of the function.
In each chamber, the index $i$ that maximizes $l_i \circ r_{N-1} \circ L_{N-1} \circ \ldots \circ r_1 \circ L_1$ is determined, and we make one outgoing arrow from the corresponding vertex to the target vertex $i\in T$.
\end{proof}	

\begin{figure}[ht]
	\includegraphics[scale=0.5]{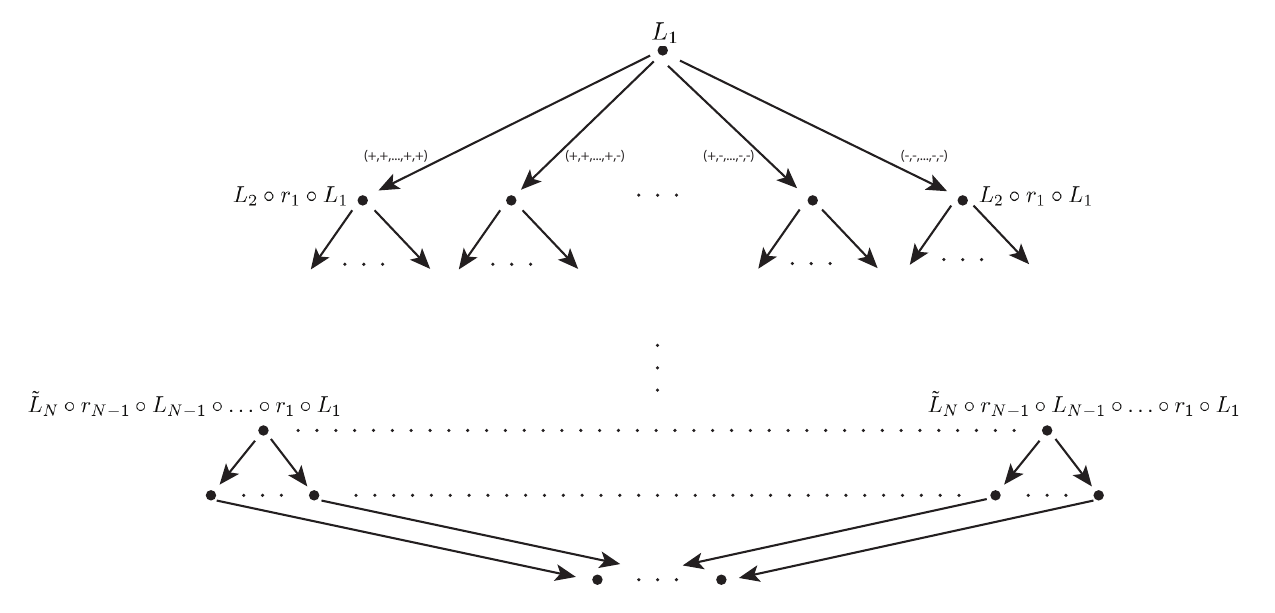}
	\caption{}
	\label{fig:linLogGrReLu}
\end{figure}

Figure \ref{fig:linLogGrReLu} depicts the logical graph in the above proposition.

In classification problems, $T$ is the set of labels for elements in $D$, and the data determines a subset in $D \times T$ as a graph of a function.  Deep learning of network models provide a way to approximate the subset as the graph of a linear logical function $\gr(f_{G,L})$.   Theoretically, this gives an interpretation of the dataset, namely, the linear logical graph gives a logical way to deduce the labels based on linear conditional statements on $D$.

The following lemma concerns about the monoidal structure on the set of linear logical functions on $D$.

\begin{lemma} \label{lem:prod}
	Let $f_{G_i,L_i}:D\to T_i$ be linear logical functions for $i \in I$ where $I=\{1,\ldots,k\}$.  Then  \[(f_{G_i,L_i}:i\in I):D \to \prod_{i\in I} T_i\] 
	is also a linear logical function.
\end{lemma}
\begin{proof}
	We construct a linear logical graph out of $G_i$ for $i\in I$ as follows.  First, take the graph $G_1$.  For each target vertex of $G_1$, we equip it with the linear function at the source vertex of $G_2$, and attach to it the graph $G_2$.  The target vertices of the resulting graph are labeled by $T_1 \times T_2$.  Similarly, each target vertex of this graph is equipped with the linear function at the source vertex of $G_3$ and attached with the graph $G_3$.
    Inductively, we obtain the required graph, whose target vertices are labeled by $\prod_{i\in I} T_i$.  By this construction, the corresponding function is $(f_{G_i,L_i}:i\in I)$.
\end{proof}

$f_{(G,L)}$ admits the following algebraic expression in the form of a sum over paths which has an important interpretation in physics.  The proof is straightforward and is omitted.  A path in a directed graph is a finite sequence of composable arrows.  The set of all linear combinations of paths and the trivial paths at vertices form an algebra by concatenation of paths.

\begin{proposition} \label{prop:pathsum}
	Given a linear logical graph $(G,L)$,
	  \[f_{(G,L)}(x) = h\left(\sum_\gamma c_\gamma(x) \,\gamma\right) = \sum_\gamma c_\gamma(x)\, h(\gamma)\] 
	where the sum is over all possible paths $\gamma$ in $G$ from the source vertex to one of the target vertices; $h(\gamma)$ denotes the target vertex that $\gamma$ heads to; for $\gamma = a_r \ldots a_1$, 
	  \[c_\gamma(x) = \prod_{i=1}^r s_{a_i}(x)\] 
	where $s_{a}(x)=1$ if $x$ lies in the chamber corresponding to the arrow $a$, or $0$ otherwise.  In the above sum, exactly one of the terms is non-zero.
\end{proposition}

\subsection{Zero locus}

Alternatively, we can formulate the graphs $\gr(f_{G,L})$ as zero loci of linear logical functions targeted at the field $\bF_2$ with two elements as follows.  Such a formulation has the advantage of making the framework of algebraic geometry available in this setting.

\begin{proposition} \label{prop:zero}
	For each linear logical function $f_{G,L}: D \to T$, there exists a linear logical function $f_{\tilde{G},\tilde{L}}: D\times T \to \bF_2$ whose zero locus in $D\times T$ equals $\gr(f_{G,L})$.
\end{proposition}

\begin{proof}
	Given a linear logical function $f_{G,L}: D \to T$, we construct another linear logical function $f_{\tilde{G},\tilde{L}}: D\times T \to \bF_2$ as follows.  Without loss of generality, let $T=\{1,\ldots,p\}$, so that $D \times T$ is embedded as a subset of $\R^{n+1}$.  Any linear function on $\R^n$ is pulled back as a linear function on $\R^{n+1}$ by the standard projection $\R^{n+1} \to \R^n$ that forgets the last component.  Then $f_{G,L}$ is lifted as a linear logical function $D\times T \to T$.  
	
	Consider the corresponding graph $(G,L)$.  For the $k$-th target vertex of $(G,L)$ (that corresponds to $k \in T$), we equip it with the linear function  \[(y-(k-1/2),(k+1/2)-y):\R^{n+1} \to \R^2\] 
	where $y$ is the last coordinate of $\R^{n+1}$.  This linear function produces three chambers in $\R^{n+1}$.  Correspondingly we make three outgoing arrows of the vertex.  Finally, the outcome vertex that corresponds to $(+,+)$ is connected to the vertex $0 \in \bF_2$; the other two outcome vertices are connected to the vertex $1 \in \bF_2$.  We obtain a linear logical graph $(\tilde{G},\tilde{L})$ and the corresponding function $f_{\tilde{G},\tilde{L}}: D\times T \to \bF_2$.
	
	By construction, $f_{\tilde{G},\tilde{L}}(x,y)=0$ for $(x,y)\in D\times T$ if and only if $y = f_{G,L}(x) \in T$.  Thus, the zero locus of $f_{\tilde{G},\tilde{L}}$ is the graph of $f_{G,L}$.
\end{proof}

The set of functions (with a fixed domain) valued in $\bF_2$ forms a unital commutative and associative algebra over $\bF_2$, which is known as a Boolean algebra.

\begin{proposition}
	The subset of linear logical functions $D \to \bF_2$ forms a Boolean ring $\cL$ (for a fixed $D \subset \R^n$).
\end{proposition}

\begin{proof}
	We need to show that the subset is closed under addition and multiplication induced from the corresponding operations of $\bF_2$.
	
	Let $f_{(G_1,L_1)}$ and $f_{(G_2,L_2)}$ be linear logical functions $D \to \bF_2$.  By Lemma \ref{lem:prod}, $(f_{(G_1,L_1)},f_{(G_2,L_2)}): D \to {(\bF_2)}^2$ is a linear logical function.  Consider the corresponding logical graph.  The target vertices are labeled by $(s_1,s_2) \in {(\bF_2)}^2$.  We connect each of them to the vertex $s_1 + s_2 \in \bF_2$ by an arrow.  This gives a linear logical graph whose corresponding function is $f_{(G_1,L_1)}+f_{(G_2,L_2)}$.  We obtain $f_{(G_1,L_1)}\cdot f_{(G_2,L_2)}$ in a similar way.
\end{proof}

In this algebro-geometric formulation, the zero locus of $f_{G,L}: D \to \bF_2$ corresponds to the ideal $(f_{G,L}) \subset \cL$.

\subsection{Parametrization}


The graph $\gr(f_{G,L})$ of a linear logical function can also be put in parametric form.  For the moment, let assume the domain $D$ to be finite.  First, we need the following lemma.

\begin{lemma} \label{lem:id}
	Assume $D\subset \R^n$ is finite.  Then the identity function $I_D: D \to D$ is a linear logical function.
\end{lemma}
\begin{proof}
Since $D$ is finite, there exists a linear function $l:\R^n \to \R^N$ such that each chamber of $l$ contains at most one point of $D$.  Then we construct the linear logical graph $G$ as follows.  The source vertex is equipped with the linear function $l$ and outgoing arrows corresponding to the chambers of $l$.  Elements of $D$ are identified as the target vertices of these arrows which correspond to chambers that contain them.  The corresponding function $f_{G,L}$ equals $I_D$.
\end{proof}

\begin{proposition}
	Given a linear logical function $f_{G,L}:D\to T$ with finite $|D|$, there exists an injective linear logical function $D \to D \times T$ whose image equals $\gr(f_{G,L})$.
\end{proposition}
\begin{proof}
	By Lemma \ref{lem:id} and \ref{lem:prod}, $(I_D, f_{G,L})$ is a linear logical function.  By definition, its image equals $\gr(f_{G,L})$.
\end{proof}

\subsection{Fuzzy linear logical functions}
Another important feature of a dataset is its fuzziness.  Below, we formulate the notion of a fuzzy linear logical function and consider its graph.  
Basic notions of fuzzy logic can be found in textbooks such as \cite{FuzzyTextbook_2013}.
There are many developed applications of fuzzy logic such as modeling, control, pattern recognition and networks, see for instance \cite{FuzzyModelingAndControl_1985, FuzzyPatternRecog_1994, FuzzyEngineeringTutorial_1995,FuzzyCNNMNIST_2016, FuzzySurvey_2020}.

\begin{definition} \label{def:fuzlinlogfcn}
	Let $G$ be a finite directed graph that has no oriented cycle, has exactly one source vertex and target vertices $t_1,\ldots,t_K$ as in Definition \ref{def:logfcn}.
	Each vertex $v$ of $G$ is equipped with a product of standard simplices
	\[P_v = \prod_{k=1}^{m_v} S^{d_{v,k}}, \,\,\, S^{d_{v,k}} = \left\{(y_0,\ldots,y_{d_{v,k}}) \in \R_{\geq 0}^{d_{v,k}+1}: \sum_{i=0}^{d_{v,k}} y_i = 1\right\}\]
	for some integers $m_v>0$, $d_{v,k}\geq 0$.  $P_v$ is called the internal state space of the vertex $v$.  Let $D$ be a subset of the internal state space of the source vertex of $G$.  Each vertex $v$ that has more than one outgoing arrows is equipped with an affine linear function 
	\[ l_v: \prod_{k=1}^{m_v} \R^{d_{v,k}} \to \R^j \]
	for some $j > 0$, whose chambers in the product simplex $P_v$ are one-to-one corresponding to the outgoing arrows of $v$.  (In above, $\R^{d_{v,k}}$ is identified with the affine subspace $\left\{\sum_{i=0}^{d_{v,k}} y_i = 1\right\}$ that contains $S^{d_{v,k}}$.)
	Let $L$ denote the collection of these affine linear functions.
	Moreover, each arrow $a$ is equipped with a continuous function 
	\[ p_a: P_{s(a)} \to P_{t(a)} \]
	where $s(a),t(a)$ denote the source and target vertices respectively.  
	
	We call $(G,L,P,p)$ a fuzzy linear logical graph.  $(G,L,P,p)$ determines a function 
	 \[f_{(G,L,P,p)}: D \to P^{\textrm{out}} := \coprod_{l=1}^K P_{t_l}\]  
	as follows.
	Given $x \in D$, as in Definition \ref{def:logfcn}, the collection $L$ of linear functions over vertices of $G$ evaluated at the image of $x$ under the arrow maps $p_a$ determines a path from the source vertex to one of the target vertices $t_l$.  By composing the corresponding arrow maps $p_a$ on the internal state spaces along the path and evaluating at $x$, we obtain a value $f_{(G,L,P,p)}(x) \in P_{t_l}$.  The resulting function $f_{(G,L,P,p)}$ is called a fuzzy linear logical function.
\end{definition}

\begin{remark} \label{rmk:fuzzy}
	Note that the linear functions $l_v$ in $L$ in the above definition have domain to be the internal state spaces over the corresponding vertices $v$.  In comparison, the linear functions $l_v$ in $L$ in Definition \ref{def:logfcn} have domain to be the input space $\R^n$.  A linear logical graph $(G,L)$ in Definition \ref{def:logfcn} has no internal state space except the $\R^n$ at the input vertex.
	
	To relate the two notions given by the above definition and Definition \ref{def:logfcn},
	we can set $P_v = S^n$ to be the same for all vertices $v$ except the target vertices $t_1,\ldots,t_K$, which are equipped with the zero-dimensional simplex (a point), and set $p_a$ to be identity maps for all arrows $a$ that are not targeted at any of $t_i$.  Then $f_{(G,L,P,p)}$ reduces back to a linear logical function in Definition \ref{def:logfcn}.
	
	We call the corners of the convex set $P_v$ to be state vertices, which takes the form $e_I = (e_{i_1},\ldots,e_{i_{m_v}}) \in P_v$ for a multi-index $I=(i_1,\ldots,i_{m_v})$, where $\{e_0,\ldots,e_{d_{v,k}}\} \subset \R^{d_{v,k}+1}$ is the standard basis.
	We can construct a bigger graph by replacing each vertex of $G$ by the collection of state vertices, and each arrow of $G$ by the collection of all possible arrows from source state vertices to target state vertices.
	Then the vertices of $G$ are interpreted as `layers' or `clusters' of vertices of this bigger graph.  The input state $x\in D$, the arrow linear functions $L$ and the maps between state spaces $p$ determine the probability of getting to each target state vertex from the source vertex.  
	
	Under this interpretation, we take the target set to be the disjoint union of corners of $P_{t_l}$ at the target vertices $t_1,\ldots,t_K$:
	\begin{equation} \label{eq:T}
		T = \coprod_{l=1}^K T_l := \coprod_{l=1}^K \left\{e_I: I = (i_1,\ldots,i_{m_{t_l}}) \textrm{ for } i_k \in \{0,\ldots,d_{t_l,k}\}\right\}
	\end{equation}
	which is a finite set.  The function $f = f_{(G,L,P,p)}$ determines the probability of the outcome for each input state $x \in D$ as follows.  Let $f(x) \in P_{t_l} = \prod_{k=1}^{m_{t_l}} S^{d_{t_l,k}}$ for some $l=1,\ldots,K$.  Then the probability of being in $T_j$ is zero for $j\not=l$.  Writing $f(x) = (f_1(x),\ldots,f_{m_{t_l}}(x))$ for $f_k(x)\in S^{d_{t_l,k}}$, the probability of the output to be $t=e_I \in T_l$ for $I=(i_1,\ldots,i_{m_{t_l}})$ is given by $\prod_{k=1}^{m_{t_l}} f^{(i_k)}_k(x) \in [0,1]$.
\end{remark}

\begin{proposition}\label{prop:sigmoid2}
	Consider the function 
	  \[f = \tilde{\sigma} \circ L_{N} \circ \tilde{s}_{N-1} \circ L_{N-1} \circ \ldots \circ \tilde{s}_1 \circ L_1. \]	
	given by a feed-forward network model whose activation function at each hidden layer is the sigmoid function (denoted by $\tilde{s}_i$), and that at the last layer is the softmax function $\tilde{\sigma}$. $f$ is a fuzzy linear logical function.
\end{proposition}

\begin{proof}
	We set $(G,L,P,p)$ as follows.  $G$ is the graph that has $(N+1)$ vertices $v_0,\ldots,v_N$ with arrows $a_i$ from $v_{i-1}$ to $v_i$ for $i=1,\ldots,N$.  $L$ is just an empty set.  $P_i := {(S^1)}^{m_i}$ for $i=0,\ldots,N-1$, where $m_i$ is the dimension of the domain of $L_{i+1}$.   The one-dimensional simplex $S^1$ is identified with the interval $[0,1]$.  $P_N := S^{m_N}$ where $m_N$ is the dimension of the target of $L_N$.  Then
\begin{align*}
	p_i := \tilde{s}_{i} \circ L_{i}|_{{[0,1]}^{m_{i-1}}} \textrm{ for } i=1,\ldots,N-1\\
	p_N := \tilde{\sigma} \circ L_{N}|_{{[0,1]}^{m_{N-1}}}.
\end{align*}
Then $f = f_{(G,L,P,p)}$.
\end{proof}
	  
\begin{proposition} \label{prop:ReLu2}
		Consider the function
	\[f = \tilde{\sigma} \circ L_{N} \circ r_{N-1} \circ L_{N-1} \circ \ldots \circ r_1 \circ L_1. \]		
		given by a feed-forward network model whose activation function at each hidden layer is the ReLu function (denoted by $r_i$), and that at the last layer is the softmax function $\tilde{\sigma}$. ($L_i$ denotes affine linear functions.) $f$ is a fuzzy linear logical function.
\end{proposition}

\begin{proof}
	We need to construct a fuzzy linear logical graph $(G,L,P,p)$ such that $f = f_{(G,L,P,p)}$.
	We take $G$ to be the logical graph constructed in Example \ref{ex:ReLu} (Figure \ref{fig:linLogGrReLu}) with the last two layers of vertices replaced by a single target vertex $t$. Each vertex that targets at the last vertex $t$ and $t$ itself only has zero or one outgoing arrow and hence is not equipped with linear function. Other vertices are equipped with linear functions on the input space $\R^n$ as in Example \ref{ex:ReLu}. 
	We take the internal state space to be the $n$-dimensional cube $P_v := {(S^1)}^{n} \subset \R^n$ where $n$ is the input dimension at every vertex $v$ except at the target vertex, whose internal state space is defined to be the simplex $P_t := S^d$ where $d$ is the target dimension of $f$.  The function $p_a$ in Definition \ref{def:fuzlinlogfcn} is defined to be the identity function on the internal state space $P_v$ for every arrow $a$ except for the arrows that target at $t$.  Now we need to define $p_a$ for the arrows that target at $t$. Let $t_i'$ be the source vertices of these arrows. The input space $\R^n$ is subdivided into chambers $\{x \in \R^n: \textrm{the path determined by } x \textrm{ targets at } t_i'\}$. Moreover, $L_N \circ r_{N-1} \circ L_{N-1} \circ \ldots \circ r_1 \circ L_1$ is a piecewise-linear function, whose restriction on each of these chambers is linear and extend to a linear function $l$ on $\R^n$.  Then $p_a$ for the corresponding arrow $a$ is defined to be $\tilde{\sigma}\circ l$.  By this construction, we have $f = f_{(G,L,P,p)}$.
\end{proof}

As in Proposition \ref{prop:pathsum}, $f_{(G,L,P,p)}$ can be expressed in the form of sum over paths.

\begin{proposition}
	Given a fuzzy linear logical graph $(G,L,P,p)$,
	  \[f_{(G,L,P,p)}(x) = \sum_\gamma c_\gamma(x)\, p_\gamma(x)\] 
	where the sum is over all possible paths $\gamma$ in $G$ from the source vertex to one of the target vertices; 
	for $\gamma = a_r \ldots a_1$, 
	$p_\gamma(x) = \prod_{i=1}^r p_{a_r}\ldots p_{a_1}(x)$;
	  \[c_\gamma(x) = \prod_{i=1}^r s_{a_i}(p_{a_{i-1}\ldots a_1}(x))\] 
	where $s_{a}(x)=1$ if $x \in P_{t_a}$ lies in the chamber corresponding to the arrow $a$, or $0$ otherwise.  In the above sum, exactly one of the terms is non-zero.
\end{proposition}

\subsection{As a fuzzy subset}
A fuzzy subset of a topological measure space $X$ is a continuous measurable function $\cF:X \to [0,1]$.  This generalizes the characteristic function of a subset.  
The interval $[0,1]$ can be equipped with the semiring structure whose addition and multiplication are taking maximum and minimum respectively.  This induces a semiring structure on the collection of fuzzy subsets $\cF:X\to [0,1]$ that plays the role of union and intersection operations.   

The graph $\gr(f)$ of a function 
 \[f: D \to \coprod_{l=1}^K P_{t_l},\] 
where $P_{t_l}$ are products of simplices as in Definition \ref{def:fuzlinlogfcn}, is defined to be the fuzzy subset in $D \times T$, where $T$ is defined by Equation \eqref{eq:T}, given by the probability at every $(x,t)\in D\times T$ determined by $f$ in Remark \ref{rmk:fuzzy}.  

The following is a fuzzy analog of Proposition \ref{prop:zero}.

\begin{proposition}
	Let $f_{(G,L,P,p)}$ be a fuzzy linear logical function.  Let $\cF: D\times T \to [0,1]$ be the characteristic function of its graph where $T$ is defined by \eqref{eq:T}.  Then $\cF$ is also a fuzzy linear logical function.
\end{proposition}
\begin{proof}
Similar to the proof of Proposition \ref{prop:zero}, we embed the finite set $T$ as the subset $\{1,\ldots,|T|\}\subset \R$.  The affine linear functions on $\R^n \supset D$ in the collection $L$ are pulled back as affine linear functions on $\R^{n+1} \supset D\times T$.  Similarly, for the input vertex $v_0$, we replace the product simplex $P_{v_0}$ by $\tilde{P}_{v_0} := P_{v_0} \times [0,|T|+1]$ (where the interval $[0,|T|+1]$ is identified with $S^1$); for the arrows $a$ tailing at $v_0$, $p_a$ are pulled back to be functions $\tilde{P}_{v_0} \to P_{h(a)}$.
Then we obtain $(\tilde{L},\tilde{P},\tilde{p})$ on $G$.

For each of the target vertices $t_l$ of $G$, we equip it with the linear function
  \[(y-3/2,y-5/2,\ldots,y-(|T|-1/2)): \R^{n+1} \to \R^{|T|-1}\]  
where $y$ is the last coordinate of the domain $\R^{n+1}$.  It divides $\R^{n+1}$ into $|T|$ chambers that contain $\R^n \times \{j\}$ for some $j=1,\ldots,|T|$.  Correspondingly we make $|T|$ outgoing arrows of the vertex $t_l$.  The new vertices are equipped with the internal state space $P_{t_l}$, and the new arrows are equipped with the identity function $P_{t_l}\to P_{t_l}$.
Then we get $|T| K$ additional vertices, where $K$ is the number of target vertices $t_l$ of $G$.  Let's label these vertices by $v_{j,l}$ for $j \in T$ and $l\in \{1,\ldots,K\}$. Each of these vertices are connected to the new output vertex by a new arrow.  The new output vertex is equipped with the internal state space $S^1\cong [0,1]$.  The arrow from $v_{j,l}$ to the output vertex is equipped with the following function $p_{j,l}: P_{v_{j,l}} \to [0,1]$.  If $j \not\in T_l$, then we set $p_{j,l}\equiv 0$.
Otherwise, for $j=e_I \in T_l$ and $I=(i_1,\ldots,i_{m_{t_l}})$, $p_{j,l} := \prod_{k=1}^{m_{t_l}} u^{(i_k)}_k$ where $u_k^{(0)},\ldots,u_k^{(d_{t_l,k})}$ are the coordinates of $S^{d_{t_l,k}} \subset \R_{\geq 0}^{d_{t_l,k}+1}$.
This gives the fuzzy linear logical graph $(\tilde{G},\tilde{L},\tilde{P},\tilde{p})$ whose associated function $f_{(\tilde{G},\tilde{L},\tilde{P},\tilde{p})}: D \times T \to [0,1]$ is the characteristic function.
\end{proof}

The above motivates us to consider fuzzy subsets whose characteristic functions are fuzzy linear logical functions $\cF:X \to [0,1]$.  Below, we will show that they form a sub-semiring, that is, they are closed under fuzzy union and intersection.  We will need the following lemma analogous to Lemma \ref{lem:prod}.

\begin{lemma} \label{lem:fuz_prod}
	Let $f_{G_i,L_i,P_i,p_i}:D\to P^{\textrm{out}}_i$ be fuzzy linear logical functions for $i \in I$ where $I=\{1,\ldots,k\}$, and assume that the input state space $P_{i,\mathrm{in}}$ are the same for all $i$.  Then  \[(f_{G_i,L_i,P_i,p_i}:i\in I):D \to \prod_{i\in I} P^{\textrm{out}}_i\] 
	is also a fuzzy linear logical function.
\end{lemma}
\begin{proof}
	By the proof of Lemma \ref{lem:prod}, we obtain a new graph $(G,L)$ from $(G_i,L_i)$ for $i=1,\ldots,k$ by attaching $(G_{i+1},L_{i+1})$ to the target vertices of $(G_i,L_i)$.
	For the internal state spaces, we change as follows.
	First, we make a new input vertex $\tilde{v}_0$ and an arrow $\tilde{a}_0$ from $\tilde{v}_0$ to the original input vertex $v_0$ of $(G,L)$.
	We denote the resulting graph by $(\tilde{G},\tilde{L})$.
	We define $\tilde{P}_{\tilde{v}_0} := P_{v_0}$, $\tilde{P}_{v_0} := \prod_{i=1}^k P_{v_0}$ where $P_{v_0} = P_{i,\mathrm{in}}$ for all $i$ by assumption, and $\tilde{p}_{\tilde{a}_0}:P_{v_0} \to \prod_{i=1}^k P_{v_0}$ to be the diagonal map $\tilde{p}_{\tilde{a}_0} = (\mathrm{Id},\ldots,\mathrm{Id})$.
	The internal state spaces $P_v$ over vertices $v$ of $G_1$ are replaced by $P_v \times \prod_{i=2}^k P_{v_0}$, and $p_a$ for arrows $a$ of $G_1$ are replaced by $(p_a,\mathrm{Id},\ldots,\mathrm{Id})$.
	Next, over vertices $v$ of the graph $G_2$ that is attached to the target vertex $t_l$ of $G_1$, the internal state space $P_v$ is replaced by $P_{t_l} \times P_v \times \prod_{i=3}^k P_{v_0}$, and $p_a$ for arrows $a$ of $G_2$ are replaced by $(\mathrm{Id},p_a,\mathrm{Id},\ldots,\mathrm{Id})$.
	Inductively, we obtain the desired graph $(\tilde{G},\tilde{L},\tilde{P},\tilde{p})$.
\end{proof}

\begin{proposition}
	Suppose $\cF_1,\cF_2:X \to [0,1]$ are fuzzy subsets defined by fuzzy linear logical functions.  Then $\cF_1\cup\cF_2$ and $\cF_1\cap\cF_2$ are also fuzzy subsets defined by fuzzy linear logical functions.
\end{proposition}

\begin{proof}
	By the previous lemma, $(\cF_1,\cF_2) = f_{(G,L,P,p)}: X \to [0,1]\times [0,1]$ for some fuzzy linear logical graph $(G,L,P,p)$, which has a single output vertex whose internal state space is $[0,1]\times [0,1] \cong S^1 \times S^1$.  We attach an arrow $a$ to this output vertex.  Over the new target vertex $v$, $P_v := [0,1]$; $p_a := \max: [0,1]\times [0,1] \to [0,1]$ (or $p_a := \min$).  Then we obtain $(\tilde{G},\tilde{L},\tilde{P},\tilde{p})$ whose corresponding fuzzy function defines $\cF_1\cup\cF_2$ (or $\cF_1\cap\cF_2$ respectively).
\end{proof}

\begin{remark} \label{rmk:fuzzy2}
	For $f:P^{\textrm{in}}\to P^{\textrm{out}}$ where $P^{\textrm{in}},P^{\textrm{out}}$ are product simplices, we can have various interpretations.  
	\begin{enumerate}
		\item As a usual function, its graph is in the product $P^{\textrm{in}}\times P^{\textrm{out}}$.
		\item As a fuzzy function on $P^{\textrm{in}}$: $P^{\textrm{in}} \to T$ where $T$ is the finite set of vertices of the product simplex $P^{\textrm{out}}$, its graph is a fuzzy subset in $P^{\textrm{in}} \times T$.
		\item The domain product simplex $P^{\textrm{in}}$ can also be understood as a collection of fuzzy points over $V$, the finite set of vertices of $P^{\textrm{in}}$, where a fuzzy point here just refers to a probability distribution (which integrates to $1$).
		\item Similarly, $P^{\textrm{in}}\times P^{\textrm{out}}$ can be understood as a collection of fuzzy points over $V\times T$.  Thus, the (usual) graph of $f$ can be interpreted as a sub-collection of fuzzy points over $V\times T$.
	\end{enumerate}
\end{remark}

$(\Id,f)$ gives a parametric description of the graph of a function $f$.  The following ensures that it is a fuzzy linear logical function if $f$ is.

\begin{cor}
	Let $f: D \to P^{\mathrm{out}}$ be a fuzzy linear logical function.  Then $(\Id,f): D \to D \times P^{\mathrm{out}}$ is also a fuzzy linear logical function whose image is the graph of $f$.
\end{cor}
\begin{proof}
	By Lemma \ref{lem:fuz_prod}, it suffices to know that $\Id: D \to D$ is a fuzzy linear logical function.  This is obvious: we take the graph with two vertices serving as input and output, which are connected by one arrow.  The input and output vertices are equipped with the internal state spaces that contain $D$, and $p$ is just defined by the identity function.
\end{proof}

\begin{remark}
	Generative deep learning models widely used nowadays can be understood as parametric descriptions of data sets $X$ by fuzzy linear logical functions $f:D\to X$ (where $D$ and $X$ are embedded in certain product simplices $P^{\mathrm{in}}$ and $P^{\mathrm{out}}$ respectively and $D$ is usually called the noise space). We focus on classification problems in the current work and plan to extend the framework to other problems as well in the future.
\end{remark}

\subsection{A digression to non-linearity in a quantum-classical system} \label{sec:quant}
The fuzzy linear logical functions in Definition \ref{def:fuzlinlogfcn} have the following quantum analog.  Quantum systems and quantum random walks are well known and studied, see for instance \cite{GeometicQuantumMechanics_1999} and \cite{QuantumRW_2003}.  On the other hand, they depend only linearly on the initial state in probability.  The motivation of this subsection is to compare fuzzy and quantum systems and to show how non-linear dependence on the initial probability distribution can come up. On the other hand, this section is mostly unrelated to the rest of the writing and can be skipped. 

\begin{definition} \label{def:quan}
	Let $G$ be a finite directed graph that has no oriented cycle and has exactly one source vertex and target vertices $t_1,\ldots,t_K$.  Each vertex $v$ is equipped with a product of projectifications of Hilbert spaces over complex numbers: 
	\[ Q_v := \prod_{l=1}^{m_v} \bP( H_{v,l}) \] 
	for some integer $m_v>0$.  We fix an orthonormal basis in each Hilbert space $H_{v,l}$, which gives a basis in the tensor product:
	 \[E^{(v)} = \left\{e^{(v)}_I = (e^{(v)}_{i_1},\ldots,e^{(v)}_{i_{m_{v}}}): e^{(v)}_{i_l} \textrm{ is a basic vector of } H_{v,l}\right\}.\] 
	For each vertex $v$ that has more than one outgoing arrows, we make a choice of a decomposition of the set $E^{(v)}$ into subsets that are one-to-one corresponding to the outgoing arrows.  Each arrow $a$ is equipped with a map $q_a$ from the corresponding subset of basic vectors $e^{(t(a))}_I$ to $Q_{h(a)}$.  
	
	Let's call the tuple $(G,Q,E,q)$ to be a quantum logical graph.  
\end{definition}

We obtain a probabilistic map $f_{(G,Q,E,q)}:Q^{\mathrm{in}} \to T$ as follows.  Given a state $\vec{w}=(w_1,\ldots,w_{m_v}) \in Q_v$ at a vertex $v$, we make a quantum measurement and $\vec{w}$ projects to one of the basic elements $e^{(v)}_I$ with probability $\prod_{l=1}^{m_v} |\langle w_l, e^{(v)}_{i_l}\rangle|^2$.  The outcome $e^{(v)}_I$ determines which outgoing arrow $a$ to pick, and the corresponding map $q_a$ sends it to an element of $Q_{h(a)}$.  Inductively we obtain an element $f_{(G,Q,E,q)}(w) \in T$.

However, such a process is simply linearly depending on the initial condition in probability: the probabilities of outcomes of the quantum process $f_{(G,Q,E,q)}(w)$ for an input state $w$ (which is complex-valued) simply linearly depends on the modulus of components of the input $w$.  In other words, the output probabilities are simply obtained by a matrix multiplication on the input probabilities.  To produce non-linear physical phenomena, we need the following extra ingredient.

Let's consider the state space $\bP^n$ of a single particle. A basis gives a map $\mu:\bP^n \to S^n$ to the simplex $S^n$ (also known as the moment map of a corresponding torus action on $\bP^n$):
\[ \mu = \left( \frac{|z_0|^2}{|z_0|^2+\ldots+|z_n|^2}, \ldots, \frac{|z_n|^2}{|z_0|^2+\ldots+|z_n|^2} \right):\bP^n \to S^n. \]
The components of the moment map are the probability of quantum projection to basic states of the particle upon observation.  By the law of large numbers, if we make independent observations of particles in an identical quantum state $\vec{z} = [z_0:\ldots:z_n] \in \bP^n$ for $N$ times, the average of the observed results (which are elements in $\{(1,0,\ldots,0),\ldots,(0,\ldots,0,1)\}$) converges to $\mu(\vec{z}) \in S^n$ as $N\to \infty$.  

The additional ingredient we need is a choice of a map $s: S^n \to \bP^m$ and $m,n \in \Z_{>0}$.  For instance, for $m=n=1$, we set the initial phase of the electron spin state according to a number in $[0,1]$. 

Upon an observation of a state, we obtain a point in $\{(1,0,\ldots,0),\ldots,(0,\ldots,0,1)\} \subset S^n$.  Now if we have $N$ particles simultaneously observed, we obtain $N$ values, whose average is again a point $p$ in the simplex $S^n$.  By $s$, these are turned to $N$ quantum particles in state $s(p) \in \bP^m$ again.

$\mu: \bP^n \to S^n$ and $s:S^n \to \bP^m$ give an interplay between quantum processes and classical processes with averaging.  Averaging in the classical world is the main ingredient to produce non-linearity from the linear quantum process. 

Now, let's modify Definition \ref{def:quan} by using $s:S^n \to \bP^m$.  Let $P_v$ be the product simplex corresponding to $Q_v$ at each vertex. 
Moreover, as in Definition \ref{def:logfcn} and \ref{def:fuzlinlogfcn} for (fuzzy) linear logical functions,  we equip each vertex with affine linear functions $l_v$ whose corresponding systems of inequalities divide $P_{v}$ into chambers.   This decomposition of $P_{v}$ plays the role of the decomposition of $E^{(v)}$ in Definition \ref{def:quan}.  The outgoing arrows at $v$ are in a one-to-one correspondence with the chambers.  Each outgoing arrow $a$ at $v$ is equipped with a map $\tilde{q}_a$ from the corresponding chamber of $P_{v}$ to $Q_{h(a)}$.  $\tilde{q}_a$ can be understood as an extension of $q_a$ (whose domain is a subset of corners of $P_{t(a)}$) in Definition \ref{def:quan}. 

\begin{definition}
	We call the tuple $(G,Q,E,L,\tilde{q})$ (where $L$ is the collection of affine linear functions $l_v$) to be a quantum-classical logical graph.
\end{definition}

Given $N$ copies of the same state $\vec{w}$ in $Q_v$, we first take a quantum projection of these and they become elements in $P_v$.  We take an average of these $N$ elements, which lies in a certain chamber defined by $l_v$.  The chamber corresponds to an outgoing arrow $a$, and the map $\tilde{q}_a$ produces $N$ elements in $Q_{h(a)}$.  Inductively we obtain a quantum-classical process $Q^{\mathrm{in}} \to T$.

For $N=1$, this essentially produces the same linear probabilistic outcomes as in Definition \ref{def:quan}.  On the other hand, when $N>1$, the process is no longer linear and produces a fuzzy linear logical function $P^{\mathrm{in}} \to P^{\mathrm{out}}$.

In summary, non-linear dependence on the initial state results from averaging of observed states.

\begin{remark}
	We can allow loops or cycles in the above definition.  Then the system may run without stop.  In this situation, the main object of concern is the resulting (possibly infinite) sequence of pairs $(v,s)$, where $v$ is a vertex of $G$ and $s$ is a state in $Q_v$.  This gives a quantum-classical walk on the graph $G$.
	
	We can make a similar generalization for (fuzzy) linear logical functions by allowing loops or cycles.  This is typical in applications in time-dependent network models.
\end{remark}

\section{Linear logical structures for a measure space}

In the previous section, we have defined linear logical functions based on a directed graph.  In this section, we will first show the equivalence between our definition of linear logical functions and semilinear functions\cite{vandenDries} in the literature.  Thus, the linear logical graph we have defined can be understood as a representation of semilinear functions.  Moreover, fuzzy and quantum logical functions that we define can be understood as deformations of semilinear functions.

Next, we will consider measurable functions and show that they can be approximated and covered by semilinear functions.  This motivates the definition of a logifold, which is a measure space that has graphs of linear logical functions as local models.

\subsection{Equivalence with semilinear functions} \label{sec:semilinear}
Let's first recall the definition of semilinear sets.

\begin{definition} \label{def:semilinear}
	For any positive integer $n$, semilinear sets are the subsets of $\R^n$ that are finite unions of sets of the form
	\begin{equation} \label{eq:semilinearSet}
		\{x\in\R^n : f_1(x) = \cdots = f_k(x) = 0, g_1(x) >0, \ldots, g_l(x)>0 \}
	\end{equation}
	where the $f_i$ and $g_j$ are affine linear functions.
	
	A function $f:D \to T$ on $D\subset \R^n$, where $T$ is a discrete set, is called to be semilinear if for every $t\in T$, $f^{-1}\{t\}$ equals to the intersection of $D$ with a semilinear set.
\end{definition}

Now let's consider linear logical functions defined in the last section.  We show that the two notions are equivalent (when the target set is finite).  Thus, a linear logical graph can be understood as a graphical representation (which is not unique) of a semi-linear function.  From this perspective, the last section provides fuzzy and quantum deformations of semi-linear functions.

\begin{theorem}\label{thm:semilinear}
	Consider $f : D \to T$ for a finite set $T=\{t_1,\ldots,t_s\}$ where $D\subset \R^n$.  $f$ is a semilinear function if and only if it is a linear logical function.
\end{theorem}
\begin{proof}
	It suffices to consider the case $D=\R^n$.
	We will use the following terminologies for convenience.
	Let $V(G)$ and $E(G)$ be the sets of vertices and arrows respectively for a directed graph $G$.
	A vertex $v \in V(G)$ is called to be \emph{nontrivial} if it has more than one outgoing arrows.
	It is said to be \emph{simple} if it has exactly one outgoing arrow.
	We call a vertex that has no outgoing arrow to be a target, and that has no incoming arrow to be a source.
	For a target $t$, let $R_t$ be the set of all paths from the source to target $t$.
	
	Consider a linear logical function $f = f_{(G,L)}: \R^n \to T$.  Let $p$ be a path in $R_t$ for some $t\in T$.  Let $\{v_1, \ldots, v_k\}$ be the set of non-trivial vertices that $p$ passes through.  This is a non-empty set unless $f$ is just a constant function (recall that $G$ has only one source vertex).
	At each of these vertices $v_i$, $\R^n$ is subdivided according to the affine linear functions $g_{i,1},\ldots,g_{i,N_i}$ into chambers $C_{i,1}, \ldots, C_{i,m_i}$ where $m_i$ is the number of its outgoing arrows.  All the chambers $C_{i,j}$ are semilinear sets.  
	
	For each path $p \in R_t$, we define a set $E_p$ such that $x\in E_p$ if $x$ follows path $p$ to get the target $t$.
	Then $E_p$ can be represented as $ C_{1,j_1}\cap \cdots \cap C_{k,j_k},$ which is semilinear. Moreover, the finite union
 	\[f^{-1}(t)=\bigcup_{p\in R_t} E_p\]
	is also a semilinear set.  This shows that $f$ is a semilinear function.
	
	Conversely, suppose that we are given a semilinear function.
	Without loss of generality, we can assume that $f$ is surjective.
	For every $t\in T$, $f^{-1}(t_i)$ is a semilinear set defined by a collection of affine linear functions in the form of \eqref{eq:semilinearSet}.
	Let $\mathcal{F} = \{l_1,\ldots,l_N\}$ be the union of these collections over all $t\in T$.
	
	Now we construct a linear logical graph associated to $f$.  We consider the chambers made by $(l_1,-l_1, \ldots, l_N,-l_N)$ by taking the intersection of the half spaces $l_i \geq 0$, $l_i < 0$, $-l_i \geq 0$, $-l_i < 0$.  We construct outgoing arrows of the source vertex associated with these chambers.  
	
	For each $t \in T$, $l_j$ occurs in defining $f^{-1}(t)$ as either one of the following ways:
	\begin{enumerate}
		\item $l_j > 0$, which is equivalent to $l_j \geq 0$ and $-l_j <0$,
		\item $l_j =0$, which is equivalent to $l_j \geq 0$ and $-l_j \geq 0$
		\item $l_j$ is not involved in defining $f^{-1}(t)$. 
	\end{enumerate}
	Thus, $f^{-1}(t)$ is a union of a sub-collection of chambers associated to the outgoing arrows.  Then we assign these outgoing arrows with the target vertex $t$.  This is well-defined since $f^{-1}(t)$ for different $t$ are disjoint to each other.  Moreover, since $\bigcup_t f^{-1}(t) = \R^n$, every outgoing arrow is associated with a certain target vertex.
	
	In summary, we have constructed a linear logical graph $G$ which produces the function $f$.
	
\end{proof}

The above equivalence between semilinear functions and linear logical functions naturally generalizes to definable functions in other types of $o$-minimal structures.
They provide the simplest class of examples in $o$-minimal structures for semi-algebraic and subanalytic geometry \cite{vandenDries}.  The topology of sub-level sets of definable functions was recently investigated in \cite{Ji-Meng-Ding}.
Let's first recall the basic definitions.

\begin{definition}\cite{vandenDries}
	A structure $\cS$ on $\R$ consists of a Boolean algebra $\cS_n$ of subsets of $\R^n$ for each $n=0,1,2,\ldots,$ such that
	\begin{enumerate}
		\item the diagonals $\{x \in \R^n : x_i = x_j\}, 1\leq i < j \leq n$ belong to $\cS_n$;
		\item $A \in \cS_m, B\in \cS_n \implies A\times B \in \cS_{m+n}$;
		\item $A \in \cS_{n+1} \implies \pi(A) \in \cS_n$, where $\pi : \R^{n+1} \to \R^n$ is the projection map defined by $\pi(x_1,\ldots,x_{n+1}) = (x_1, \ldots, x_n)$;
		\item the ordering $\{(x,y)\in\R^2 : x<y\}$ of $\R$ belongs to $\cS_2$.
	\end{enumerate}
	A structure $\cS$ is $o$-minimal if the sets in $\cS_1$ are exactly the subsets of $\R$ that have only finitely many connected components, that is, the finite unions of intervals and points.
	
	Given a collection $\cA$ of subsets of the Cartesian spaces $\R^n$ for various $n$, such that the ordering $\{(x,y) : x<y\}$ belongs to $\cA$, define $\Def(\cA)$ as the smallest structure on the real line containing $\cA$ by adding the diagonals to $\cA$ and closing off under Boolean operations, cartesian products, and projections.
	Sets in $\Def(\cA)$ are said to be definable from $\cA$ or simply definable if $\cA$ is clear from context.
	
	Given definable sets $A \subset \R^m$ and $B\subset \R^n$ we say that a map $f : A \to B$ is definable if its graph $\Gamma(f) = \{(x,f(x)) \in \R^{m+n} : x\in A\}$ is definable.

\end{definition}

\begin{remark}
If $\cA$ consists of the ordering, the singletons $\{r\}$ for any $r\in\R$, the graph in $\R^2$ of scalar multiplications maps $x \mapsto \lambda x : \R \to \R$ for any $\lambda\in\R$, and the graph of addition $\{(x,y,z)\in\R^3: z=x+y\}$.
Then $\Def(\cA)$ consists of semilinear sets for various positive integers $n$ (Definition \ref{def:semilinear}).

Similarly, if $\cA$ consist of the ordering, singletons, and the graphs of addition and multiplication, then $\Def(\cA)$ consists of semi-algebraic sets, which are finite unions of sets of the form
\[ \left\{ x \in \R^n : f(x) = 0, g_1(x)>0, \ldots, g_l(x)>0 \right\}\]
where $f$ and $g_1,\ldots,g_l$ are real polynomials in $n$ variables, due to the Tarski-Seidenberg Theorem\cite{analyticBEMP}.

One obtains semi-analytic sets in which the above $f,g_1,\ldots,g_l$ become real analytic functions by including graphs of analytic functions.
Let \textbf{an} be the collection $\cA$ and of the functions $f : \R^n \to \R$ for all positive integers $n$ such that $f|_{I^n}$ is analytic, $I = [-1,1] \subset \R$, and $f$ is identically $0$ outside the cubes.
The theory of semi-analytic sets and subanalytic sets show that $\Def(\textbf{an})$ is $o$-minimal, and relatively compact semi-analytic sets have only finitely many connected components.
See \cite{analyticBEMP} for efficient exposition of the \L ojasiewicz-Gabrielov-Hironaka theory of semi- and subanalytic sets.
\end{remark}

\begin{theorem}
	Let's replace the collection of affine linear functions at vertices in Definition \ref{def:logfcn} by polynomials and call the resulting functions to be polynomial logical functions.
	Then $f:D \to T$ for a finite set $T$ is a semi-algebraic function if and only if $f$ is polynomial logical functions.
\end{theorem}

The proof of the above theorem is similar to that of Theorem \ref{thm:semilinear} and hence omitted.
\subsection{Approximation of measurable functions by linear logical functions} \label{sec:universal}
We consider measurable functions $f: D \to T$, where $0<\mu\left(D\right)<\infty$ and $T$ is a finite set.
The following approximation theorem for measurable functions has two distinct features since $T$ is a finite set.
First, the functions under consideration, and linear logical functions that we use, are discontinuous.
Second, the `approximating function' actually exactly equals to the target function in a large part of $D$.
Compared to traditional approximation methods, linear logical functions have an advantage of being representable by logical graphs, which have fuzzy or quantum generalizations.

\begin{theorem}[Universal approximation theorem for measurable functions]\label{thm:UAT}
	Let $\mu$ be the standard Lebesgue measure on $\R^n$.  Let $f:D\to T$ be a measurable function with $\mu(D)<\infty$ and a finite target set $T$.
	For any $\epsilon>0$, there exists a linear logical function $L: D \to \widetilde{T}$, where $\widetilde{T} = T\cup \{\ast\}$ is $T$ adjunct with a singleton, and a measurable set $E \subset D$ with $\mu(E) < \epsilon$ such that $L|_{D-E} \equiv f|_{D-E}$.

\end{theorem}
\begin{proof}
	Let $\mathcal{R} := \left\{\prod_{k=1}^{n}(a_k,b_k] \subset \R^n: a_k < b_k \textrm{ for all k} \right\}$ be the family of rectangles in $\R^n$.
	We will use the well-known fact that for any measurable set $U$ of finite Lebesgue measure, there exists a finite subcollection $\left\{R_j: j =1,\ldots,N \right\}$ of $\mathcal{R}$ such that $\mu\left(U \,\triangle\, \bigcup_1^N R_j \right) < \epsilon$ (see for instance \cite{follandRA}).
	Here, $A \,\triangle\, B := (A-B) \cup (B-A)$ denotes the symmetric difference of two subsets $A,B$.
	
	Suppose that a measurable function $f : D \to T=\{ t_1, \ldots, t_m \}$ and $\epsilon>0$ be given.
	For each $t\in T$, let $\cS_t$ be a union of finitely many rectangles of $\mathcal{R}$ that approximates $f^{-1}(t) \subset D$ (that has finite measure) in the sense that $\mu\left((f^{-1}(t) \,\triangle\, \cS_t) \right) < \frac{\epsilon}{|T|}$.
	Note that $\cS_t$ is a semilinear set.
	
	The case $m=1$ is trivial.
	Suppose $m>1$.
	Define semilinear sets $\cS_{\ast} = \underset{i<j}{\bigcup}\left(\cS_{t_i}\cap \cS_{t_j}\right)$ and $\mathfrak{S}_t = \cS_t \setminus \cS_{\ast}$ for each $t\in T$.
	Now we define $L: D \to \widetilde{T}$,
	\[
	L(p) = 
	\begin{cases}
		t & \text{if } p \in \mathfrak{S}_t \cap D\\
		\ast & \text{if } p \in D \setminus \underset{t\in T}{\bigcup} \mathfrak{S}_t
	\end{cases}\]
	which is a semilinear function on $D$.
	
	If $p \in \cS_\ast \cap D$ then $p \in \cS_{t_i}\cap \cS_{t_j}$ for some $t_i, t_j \in T$ with $t_i \neq t_j$.  In such a case, $p\in f^{-1}(t_i)$ implies $p \notin f^{-1}(t_j)$.
	It shows $ \cS_\ast \cap D \subset \underset{t\in T}{\bigcup}\left( \cS_t \setminus f^{-1}(t)\right)$.
	Also we have $D \setminus \underset{t\in T}{\bigcup} \cS_t \subset \underset{t\in T}{\bigcup} \left(f^{-1}(t)\setminus \cS_t\right)$.
	Therefore 
	\[ D \setminus \underset{t\in T}{\bigcup} \mathfrak{S}_t = \left(D \setminus  \underset{t\in T}{\bigcup}  \cS_t \right) \cup \left(D \cap \cS_\ast \right)  \subset  \underset{t\in T}{\bigcup} f^{-1}(t)\,\triangle\, \cS_t\]
	and hence 
	\[\mu\left(D \setminus \underset{t\in T}{\bigcup} \mathfrak{S}_t\right) \leq \sum_{t\in T} \mu\left(f^{-1}(t) \,\triangle\, \cS_t\right)<\epsilon. \]
	By Theorem \ref{thm:semilinear}, $L$ is a linear logical function.
\end{proof}
	
\begin{cor} \label{cor:cover}
	Let $f : D \to T$ be a measurable function where $D \subset \R^n$ is of finite measure and $T$ is finite.
	Then there exists a family $\mathcal{L}$ of linear logical functions $L_i: D_i \to T$, where $D_i \subset D$ and $L_i \equiv f|_{D_i}$, such that $ D\setminus  \underset{i}{\bigcup}\,  D_i $ is measure zero set.
\end{cor}

	
\subsection{Linear logifold} \label{sec:logifold}

To be more flexible, we can work with Hausdorff measure which is recalled as follows.

\begin{definition}
 	Let $p\geq 0$, $\delta>0$.  For any $U\subset \R^n$, $\Diam(U)$ denotes the diameter of $U$ defined by the supremum of distance of any two points in $U$.  For a subset 
	$E \subset \R^n$,
	define
	$$H^p_\delta\left(E\right) := \underset{\cU_\delta(E)}{\inf}\, \sum_{U\in \cU_\delta(E)} \Diam\left(U\right)^p  $$
	where $\cU_\delta(E)$ denotes a cover of $E$ by sets $U$ with $\Diam(U)<\delta$.
	Then the $p$-dimensional Hausdorff measure is defined as
	$H^p(E):= \lim\limits_{\delta \to 0}H^p_\delta(E)$.  The Hausdorff dimension of $E$ is
	$\dim_H(E) := \underset{p}{\inf} \{p \in [0,\infty) : H^p(E) =0\}$.
\end{definition}


\begin{definition}\label{def:linlogifold} 
	A linear logifold is a pair $(X,\cU)$, where $X$ is a topological space equipped with a $\sigma$-algebra and a measure $\mu$, $\cU$ is a collection of pairs $(U_i,\phi_i)$ where $U_i$ are subsets of $X$ such that 
	$\mu(U_i) > 0$
	and $\mu(X-\bigcup_i U_i)=0$; $\phi_i$ are measure-preserving homeomorphisms between $U_i$ and the graphs of linear logical functions $f_i: D_i \to T_i$ (with an induced Hausdorff measure), where $D_i \subset \R^{n_i}$ are $H^{p_i}$-measurable subsets in certain dimension $p_i$, and $T_i$ are discrete sets.
	
	The elements of $\cU$ are called charts.  A chart $(U,\phi)$ is called to be entire up to measure $\epsilon$ if $\mu(X-U)<\epsilon$.
\end{definition}

Comparing to a topological manifold, we require $\mu(U_i) > 0$  in place of openness condition.
Local models are now taken to be graphs of linear logical functions in place of open subsets of Euclidean spaces.

Then the results in the last subsection can be rephrased as follows.

\begin{cor}
	Let $f: D \to T$ be a measurable function on a measurable set $D \subset \R^n$ of finite measure with a finite target set $T$.  For any $\epsilon > 0$, its graph $\gr(f) \subset D \times T$ can be equipped with  a linear logifold structure that has an entire chart up to measure $\epsilon$.
\end{cor}

\begin{remark}
	In\cite{AJ}, relations between neural networks and quiver representations were studied.  In \cite{JL,JL2}, a network model is formulated as a framed quiver representation; 
	learning of the model was formulated as a stochastic gradient descent over the corresponding moduli space.
	In this language, we now take several quivers, and we glue their representations together (in a non-linear way) to form a `logifold'.  
\end{remark}

In a similar manner, we define a fuzzy linear logifold below.  By Remark \ref{rmk:fuzzy} and (2) of Remark \ref{rmk:fuzzy2}, a fuzzy linear logical function has a graph as a fuzzy subset of $D\times T$.  We are going to use the fuzzy graph as a local model for a fuzzy space $(X,\cP)$.  

\begin{definition} \label{def:fuzlogifold}
	A fuzzy linear logifold is a tuple $(X,\cP,\cU)$, where 
	\begin{enumerate}
		\item $X$ is a topological space equipped with a measure $\mu$; 
		\item $\cP: X \to (0,1]$ is a continuous measurable function;
		\item $\cU$ is a collection of tuples $(\rho_i,\phi_i, f_i)$, where $\rho_i$ are measurable functions $\rho_i: X \to [0,1]$ with $\sum_i \rho_i \leq 1_X$ that describe fuzzy subsets of $X$, 
		whose supports are denoted by by $U_i = \{x\in X: \rho_i(x)>0\} \subset X$; $$\phi_i: U_i \to D_i \times T_i$$ are measure-preserving homeomorphisms where $T_i$ are finite sets in the form of \eqref{eq:T} and $D_i \subset \R^{n_i}$ are $H^{p_i}$-measurable subsets in certain dimension $p_i$; $f_i$ are fuzzy linear logical functions on $D_i$ whose target sets are $T_i$ described in Remark \ref{rmk:fuzzy};
		\item the induced fuzzy graphs $\cF_i: D_i \times T_i \to [0,1]$ of $f_i$ satisfy
		\begin{equation} \label{eq:P}
			\cP = \sum_i \rho_i\cdot \phi_i^*(\cF_i). 
		\end{equation}
	\end{enumerate} 
\end{definition}


Persistent homology \cite{ Edelsbrunner-Letscher-Zomorodian, Zomorodian-Carlsson, Carlsson-Zomorodian-Collins-Guibas} can be defined for fuzzy spaces $(X,\cP)$ by using the filtration $X_c := \{x \in X: \cP(x) \geq c \} \subset X$ for $c\in [0,1]$ associated to $\cP$. We plan to study persistent homology for fuzzy logifolds in a future work.


	

\section{Ensemble learning and logifolds}\label{sec:algo}

In this section, we briefly review ensemble learning method\cite{EnsembleIntro} and make a mathematical formulation via logifolds. Moreover, in Section \ref{sec:notations} and \ref{section:voting}, we introduce the concept of fuzzy domain and develop a refined voting method based on this.
We view each trained model as a chart given by a fuzzy linear logical function.
The domain of each model can be a proper subset of its feature space defined by the inverse image of a proper subset of the target classes.
For each trained model, a fuzzy domain is defined using the certainty score for each input, and only inputs which lie in its certain part are accepted. In \cite{Jung-Lau}, we demonstrated in experiments that this method produces improvements in accuracy compared to taking average of outputs.

\subsection{Mathematical Description of Neural Network Learning}
Consider a subset $\mathcal{Z}$ of  $\R^n \times T$ where $T = \{c_1, \ldots, c_N\}$, which we take as the domain of a model. 
$\R^n$ is typically referred to as the \emph{feature space}, while each $c_i \in T$ represents a class. We embed $T$ as the corners of the standard simplex $S^{N-1} \subset \R^N$. One wants to find an expression of the probability distribution of $\mathcal{Z}$ in terms of a function produced by a neural network.

\begin{definition}\label{def:neural_net_function}
	The underlying graph of a neural network is a finite directed graph $G$. Each vertex $v$ is associated with $\R^{n_v}$ for some $n_v \in \Z_{>0}$, together with a non-linear function $\R^{n_v} \to \R^{n_v}$ called an activation function. 
	
	Let $\Theta$ be the vector space of linear representations. A linear representation associates each arrow $a$ with a linear map $\R^{n_{s(a)}} \to \R^{n_{t(a)}}$, where $s(a),t(a)$ are the source and target vertices respectively.
	
	Let's fix $\gamma$ to be a linear combination of paths between two fixed vertices $s$ and $t$ in $G$. The associated network function $f_{\theta}: \R^{n_s} \to \R^{n_t}$ for each $\theta \in \Theta$ is defined to be the corresponding function obtained by the sum of compositions of linear functions and activation functions along the paths of $\gamma$.
\end{definition}

One would like to minimize the function $C_\mathcal{Z} : \Theta \to \R$, $$C_\mathcal{Z}(\theta) :=  \sum_{(x,y)\in \mathcal{Z}} \|f_\theta(x) - y\|^2$$
which measures the distance between the graph of $f_\theta$ and $\mathcal{Z}$.
To do this, one takes a stochastic gradient descent(SGD) over $\Theta$. In a discrete setting, it is given by the following equation:
\[
\theta_{k+1} = \theta_k - \eta \nabla C_\mathcal{Z} -\eta W_k
\]
where $\eta \in R_{>0}$ is called \emph{step size} or \emph{learning rate}, $W_k$ is the $k^\textrm{th}$ \emph{noise} or \emph{Brownian Motion}, and $\nabla C_\mathcal{Z}$ denotes the gradient vector field of $C_\mathcal{Z}$. (In practice, the sample $\mathcal{Z}$ is divided into batches and $C$ is a sum over a batch.)

For practical purpose, the completion of the computational process is marked by the verification of \emph{epochs}.
Then the hyper-parameter space $\mathcal{H}$ for SGD is a subspace $\{(\eta, \textrm{Batch size}, \textrm{Epochs}, \textrm{Noise})\} \subset \R^3 \times \mathcal{D}$, where $\mathcal{D}$ is the space of $\R$-valued random variables with zero mean and finite variance.
This process is called training procedure, and the resulting function $g = f_{\theta_*}$ is called a trained model where $\theta_*$ be the minimizer. The $\arg\max g(x)$ is called the prediction of $g$ at $x$, which is well-defined almost everywhere. For $c_i \in T$, $g(x,c_i):= g(x)_i$ is called the certainty of the model for $x$ being in class $c_i$.

\subsection{A brief description of Ensemble Machine Learning}

Ensemble machine learning utilizes more than one classifiers to make decisions.
Dasarathy and Sheela\cite{MultipleClassifiers1979} were early contributors to this theory, who proposed the partitioning feature space using multiple classifiers.
Ensemble systems offer several advantages, including smoothing decision boundaries, reducing classifier bias, and addressing issues related to data volume.
A widely accepted key to a successful ensemble system is achieving diversity among its classifiers.
\cite{EnsembleML, Ensemble_Survey_DYCSM, Ensemble_Survey_YYHN} provide good reviews of this theory.

Broadly speaking, designing an ensemble system involves determining how to obtain classifiers with diversity and how to combine their predictions effectively.
Here, we briefly introduce popular methods.
Bagging\cite{Bagging}, short for \emph{bootstrap aggregating} trains multiple classifiers, each on a randomly sampled subset of the training dataset.
Boosting, such as AdaBoost(Adaptive Boosting)\cite{Boosting, AdaBoost} itertatively trains classifiers by focusing on the instances they misclassified in previous rounds.
In Mixture of Experts\cite{AdaptiveMixtureofLocalExperts}, each classifier specializes in different tasks or subsets of dataset, with a `gating' layer, which determines weights for the combination of classifiers.

Given multiple classifiers, an ensemble system makes decisions based on predictions from diverse classifiers and a rule for combining predictions is necessary. This is usually done by taking a weighted sum of the predictions, see for instance \cite{KUNCHEVA_majority_voting}. Moreover, the weights may also be tuned via a training process.

\subsection{Logifold structure}

Let $(X,\cP)$ be a fuzzy topological measure space with $\mu(X)>0$ where $\mu$ is the measure of $X$, which is taken as an idealistic dataset. For instance, it can be the set of all possible appearances of cats, dogs and eggs. A sample fuzzy subset $U$ of $X$ is taken and is identified with a subset of  $\R^{n} \times T$. This identification is denoted by $\phi: U \to \R^n \times T$ and $\mathcal{Z} := \phi(U)$. For instance, this can be obtained by taking pictures in a certain number of pixels for some cats and dogs, and $T$ is taken as the subset of labels $\{$`C',`D'$\}$. By the mathematical procedure given above, we obtain a trained model $g$, which is a fuzzy linear logical function, denoted by $g = g_{(G,L,P,p)}: \R^n \to S^{|T|-1}$, where $G$ is the neural network with one target vertex, $L$ and $p$ are the affine linear maps and activation functions respectively.
This is the concept of a chart of $X$ in Definition \ref{def:fuzlogifold}.

Let $g_i : \R^{n_i} \to S^{|T_i|-1}$ be a number of trained models and $G_i : \R^{n_i} \times T_i \to [0,1]$ be the corresponding certainty functions. Note that $n_i$ and $T_i$ can be distinct for different models. Their results are combined according to certain weight functions $\rho_i: X \to [0,1]$ and we get
\[ \sum_{i} \rho_i(x)G_i(\phi_i(x)): X \to [0,1]\]
where the sum is over those $i$ whose corresponding charts $U_i$ contain $x$.
This gives $\cP$ in \eqref{eq:P}, and we obtain a fuzzy linear logifold.

In the following sections, we introduce the implementation detail for finding fuzzy domains and the corresponding voting system for models with different domains.

\subsection{Thick targets and specialization}\label{sec:notations}
We consider fuzzy subsets $U_i\subset X$ and $\mathcal{Z}_i \subset \R^{n_i} \times T_i$ with identification $\phi_i : U_i \to \cZ_i$. A common fuzziness that we make use of comes from `thick targets'. For instance, let $T=\{`C\text{'},`D\text{'},`E\text{'}\}$ (continuing the example used in the last subsection). Consider $\widetilde{T} = \{\{`C\text{'},`D\text{'}\},\{`E\text{'}\}\}$, which consists of the two classes $\{`C\text{'}, `D\text{'}\}$ and $\{`E\text{'}\}$. We take a sample $\widetilde{\cZ} \subset \R^n \times \widetilde{T}$ consisting of pictures of cats, dogs and eggs with the two possible labels `cats or dogs' and `eggs'. Then we train a model with two targets ($|\widetilde{T}|=2$), and obtain $\widetilde{G}: \R^n \times \widetilde{T} \to [0,1]$.
\begin{definition}
	Let $T$ be a finite set and $\widetilde{T}$ be a subset of the power set $\mathcal{P}\left(T\right)$ such that $\emptyset\not\in \widetilde{T}$ and $t \cap t^\prime=\emptyset$ for all distinct $t, t^\prime \in \widetilde{T}$.
	
	\begin{enumerate}
		\item $t \in \widetilde{T}$ is called thin (or fine) if it is a singleton and thick otherwise.
		\item The union $\bigcup \widetilde{T} \subset T$ is called the flattening of $\widetilde{T}$.
		\item We say that $\widetilde{T}$ is full if its flattening is $T$, and fine if all its elements are thin.
	\end{enumerate}
\end{definition}

Given a model $g_i: \R^{n_i} \to S^{|\widetilde{T}_i|-1}$, where $g_i =  (g_{i,1}, \ldots, g_{i,|\widetilde{T}_i|})$, define the certainty function of $g_i$ as $C_i := \max_j g_{i,j}$. For $\alpha \in [0,1]$,
$$\cZ_{\alpha,i} := \{(x,y) \in \cZ_i: C_i(x) \geq \alpha \} $$ 
is called \emph{the certain part of  the model $g_i$}, or \emph{fuzzy domain of $g_i$}, at certainty threshold $\alpha$ in $\cZ_i$.
Let $\cM$ denote the collection of trained model $\{g_i\}_{i\in \cI}$ with identifications $\phi_i : U_i \subset X \to \cZ_i$. The union $X_\alpha = \bigcup_{i \in \mathcal{I}}\phi_i^{-1}(\mathcal{Z}_{\alpha,i})$ is called \emph{the certain part with certainty threshold $\alpha$ of $\cM$}. For instance, in the dataset of appearances of dogs, cats and eggs, suppose we have a model $g_i$ with target $T_i = \{`C\text{'},`D\text{'}\}$. $\phi_i^{-1}(\mathcal{Z}_{\alpha,i})$ is the subset of appearance of cats and dogs sampled by the set of labeled pictures $\mathcal{Z}_i \subset \R^{n_i}\times T_i$ that has certainty $\geq \alpha$ by the model. Note that as $\alpha$ decreases, there must be more or equal number of models satisfying the conditions; in particular $X_\alpha \subset X_{\alpha'}$ for $\alpha' < \alpha$. 

Table \ref{tab:symbolsforalgorithm} summarizes the notations introduced here.
\begin{table}[h]
	\centering
	\begin{tabular}{c|c}
		\toprule
		$X$ & A fuzzy topological measure space (dataset)\\		
		\midrule[0.01mm]
		$\mathcal{M}$ & A collection of trained models $\{g_i\}_{i \in \mathcal{I}}$ \\
		\midrule[0.001mm]
		$(U_i, \phi_i, g_i)$ & The identification $\phi_i : U_i \subset X \to \R^{n_i} \times \widetilde{T}_i$, \\
		& $\mathcal{Z}_i := \phi_i(U_i)$ and $\widetilde{T}_i$ the target of $g_i$.\\
		\midrule[0.001mm]
		$T_i := \bigcup \widetilde{T}_i$ &  Flattening of targets of the $i$-th model\\
		\midrule[0.001mm]
		$g_{i,j} := j-\textrm{th component of } g_i$ & The certainty of $g_i$ to $t_j \in \widetilde{T}_i$ \\
		\midrule[0.001mm]
		$C_i := \max_j g_{i,j},$ & The certainty function of $g_i$\\
		\midrule[0.001mm]
		$P_i := t_{\arg\max_j g_{i,j}}$ & The prediction function of $g_i$\\
		\midrule[0.001mm]
		$\mathcal{Z}_{\alpha,i} \subset \R^{n_i} \times T_i$ & The certain part by $g_i$ \\ &at certainty threshold $\alpha$ in $\cZ_i$\\
		
		\midrule[0.001mm]
		$X_{\alpha}= \bigcup_{i\in \cI} \phi_i^{-1}(\mathcal{Z}_{\alpha,i}) \subset X$ & The certain part of $X$ \\ &at certainty threshold $\alpha$ by $\cM$\\
		\bottomrule
	\end{tabular}
	\caption{Table of frequently used symbols}\label{tab:symbolsforalgorithm}
\end{table}

One effective method we use to generate more charts of different types is to change the target of a model, which we call \emph{specialization}\label{specialization}.
A network function $f_\theta : \R^n \to \R^{|\widetilde{T}|}$ can be turned into a `specialist' for a new target classes $\widetilde{T}' = \{t_1,\ldots,t_m\}$ where each new target $t_i \in \widetilde{T}'$ is a proper subset of target classes $\widetilde{T}$ of $f$ such that $t_i \cap t_j = \emptyset$ if $t_i \neq t_j$.

Let $G$ and $t$ be the underlying graph of a neural network function and associated target vertex.
By adding one more vertex $u$ and adjoining it to the target vertex $t$ of $G$, we can associate a function which is the composition of linear and activation functions along the arrow $a$ whose the source is $t$ and the target is $u$, with $u$ associated to $\R^{m}$. This results in a network function $g_{\theta'} : \R^{|\widetilde{T}|} \to \R^{m}$ with underlying graph $t \overset{a}{\longrightarrow} u$.
By composing $f$ and $g$, we obtain $\tilde{f}_{(\theta',\theta)} =g_{\theta'} \circ f_\theta $, whose the target classes are $t_1,\ldots,t_m$, with the concatenated graph consisting of $G$ and  $t \overset{a}{\longrightarrow} u$. Training the obtained network function $\tilde{f}_{(\theta',\theta)}$ is called a \emph{specialization}.


\subsection{Voting system}\label{section:voting}

We assume the above setup and notations for a dataset $X$ and a collection of trained models $\cM$. We introduce a voting system that utilizes fuzzy domains and incorporates trained models with different targets. 

Let $T = \{ c_1, \ldots, c_N\}$ be a given set of target classes, and suppose a measurable function $f: X \to T$ is given. In practice, this function is inferred from a statistical sample.
We will compare $f$ with the prediction obtained from $\cM$.

Consider the subcollection of models $g_i$ in $\cM$ which have flattened targets $T_i \subset T$. By abuse of notation, we still denote this subcollection by $\cM$.
We assume the following.
\begin{enumerate}
	\item If the flattened target $T_i$ of a model is minimal in the sense that $T_j \not\subset T_i$ for any other $j \not=i$, then $\tilde{T}_i$ is fine, that is, all its elements are singleton. 
	\item Every target classes $t_1,\ldots, t_k$ in a target set $\{t_1,\ldots, t_k\}$ has no intersection.
	\item $\cM$ has a trained model whose flattened target equals $T = \{c_1,\ldots, c_N\}$.
\end{enumerate}

Below, we first define a target graph, in which each node corresponds to a flattened target $T'$. Next, we consider predictions from the collection of models with the same flattened target $T_i = T'$ for some $T'\subset T$. We then combine predictions from nodes along a path in a directed graph with no oriented cycle.

\subsubsection{Construction of the target graph.}\label{implementation:target graph}

We assign a partial order to the collection of trained models $\cM$, where the weighted answers from each trained model accumulate according to this order.
The partial order is encoded by a directed graph with no oriented cycle, which we call a \emph{target graph}.
Define $\mathcal{T}$ as the collection of flattenings, that is
\[\mathcal{T} := \left\{T_i : i \in \mathcal{I}	\right\}.\]
Among the flattenings in $\mathcal{T}$, the partially ordered subset relation induces a directed graph that has a single source vertex, called the \emph{root node}, which is associated with the given set of target classes $T = \{c_1, \ldots, c_N\}$.


Let $\mathfrak{T}$ denote the set of nodes in the target graph.
For each node $s \in \mathfrak{T}$, let $T_s \in \mathcal{T}$ denote the associated flattening of the target, and define $\mathcal{I}_s$ as the index set
\[
\mathcal{I}_s := \{ i \in \mathcal{I} \mid T_i = T_s \},
\]
which records the indices of the trained models corresponding to node $s$.
By abuse of notation, $\mathfrak{T}$ also refers to the target graph itself, and $g_i \in \mathcal{I}_v$ indicates that $g_i$ is the trained model whose index $i$ belongs to $\mathcal{I}_v$.

We define the following \emph{the refinement of targets}.
Refinements allow us to systematically combine the predictions from multiple models at each node.

\begin{definition}
	Let $T$ be a finite set and suppose we have a collection of subsets $\widetilde{T}_i$ of the power set $\cP(T)$ for $i \in \cI$ such that for each $i$, its flattening equals $T$, that is $\bigcup \widetilde{T}_i = T$; moreover $\emptyset\not\in \widetilde{T}_i$ and $t \cap t^\prime=\emptyset$ for all distinct $t, t^\prime \in \widetilde{T}_i$. The common refinement is defined to be
	$$ \left\{\bigcap_{i\in \cI} t_i: (t_i)_{i \in \cI} \in \prod_{i \in \cI} \widetilde{T}_i \right\} \setminus \left\{\emptyset\right\}.$$
\end{definition}

At each node $v \in \mathfrak{T}$, we consider the collection $\{\widetilde{T}_i\}_{i \in \mathcal{I}_v}$ of targets of all models at the node, and take its 
common refinement $\overline{T}_v$. See Example \ref{ex:refinement}.


\subsubsection{Voting rule for multiple models sharing the same target}
Let $I_E$ denote the characteristic function of a measurable set $E$, defined as
\[
I_E(x) = 
\begin{cases*}
	1  \quad \textrm{ if } x \in E, \\
	0  \quad \textrm{ otherwise},
\end{cases*}
\]
where $E \subset X$ or $E \subset \R^n$ for some positive integer $n$.

Let $(U_i,\phi_i,g_i)$ be a triple consisting of trained model $g_i : \R^{n_i} \to \widetilde{T}_i$, fuzzy subset $U_i \subset X$, and identification $\phi_i : U_i \to \cZ_i \subset \R^{n_i}\times \widetilde{T}_i$.
Let 
\begin{align*}
	\hat{x}_i := & \pi_1 \circ \phi_i(x) \in \R^{n_i}, \\ \hat{y}_i := & \pi_2 \circ \phi_i(x) \in \widetilde{T}_i
\end{align*}
denote the feature and output of realized data for each $x \in U_i$ via identification $\phi_i$ with the projection maps $\pi_1$ and $\pi_2$ from $ \R^{n_i} \times \widetilde{T}_i $ onto their first and second components, respectively.

\emph{The accuracy function $\Phi_i$ of the trained model $g_i$ over $\mathcal{Z}_{i}$ at certainty threshold $\alpha$} is defined as
\[
\Phi_{i}(\alpha) := \frac{ \left| \{ (\hat{x}_i,\hat{y}_i) \in \mathcal{Z}_{\alpha,i} : P_i(\hat{x}_i) = \hat{y}_i\}\right| }{\left| \mathcal{Z}_{\alpha,i}\right| }.
\]
Here, we denote the measure of a subset $Z$ by $|Z|$.

Let $\mathcal{G} = \{g_1, \ldots, g_m\}$ be a family of trained models sharing the same target set $\widetilde{T} = \{t_1,\ldots, t_k\}$ with accuracies $\Phi_1, \ldots, \Phi_m$, respectively. We define \emph{the weighted answer from $\mathcal{G}$, a group of trained models sharing the same targets, for $x \in X$ at certainty threshold $\alpha$} as 
\[
	\bold{p}_{\mathcal{G}}(\alpha,x) := \left(\bold{p}_{\mathcal{G}}(\alpha,x)_1, \ldots, \bold{p}_{\mathcal{G}}(\alpha,x)_k\right) := 
	\left(\frac{\sum_i I_{X_{\alpha,i}}(x) \Phi_i(\alpha) g_{i,j}(\hat{x}_i)}{\Phi_\mathcal{G}(\alpha)}\right)_{j =1 ,\ldots, k},
\]
where  $\Phi_\mathcal{G}(\alpha) := \sum_i I_{\{C_i(p_i) \geq \alpha\}}(p_i) \Phi_i(\alpha)$ is \emph{accuracy functions of models in $\mathcal{G}$} with certainty threshold $\alpha$, and $X_{\alpha,i}:= \phi^{-1}\left(\cZ_{\alpha,i}\right)$ is the certain part of $g_i$ in $X$ with certainty threshold $\alpha$ for each $i = 1,\ldots, m$. If $\Phi_\mathcal{G}( \alpha) = 0$, then we define $\bold{p}_{\mathcal{G}}(\alpha,x)  = \mathbf{0} \in \R^k$. See Example~\ref{ex:refinement}.
%
%

\subsubsection{Voting rule at a node}

Let $v$ be a node in the target graph $\mathfrak{T}$ associated with flattened target $T_v$, refinements $\overline{T}_v$, and associated models $\mathcal{I}_v$. Consider the collection of all distinct target sets of models in $ \mathcal{I}_v$, that is $\{\widetilde{T}_i\}_{i \in \mathcal{I}_v} = \{\widetilde{T}_{v,1}, \ldots, \widetilde{T}_{v,n_v}\}$ for some positive integer $n_v$. Let $\mathcal{G}_{i}$ denote the family of models sharing the same target $\widetilde{T}_{v,i} = \{t_{i,1}, \ldots, t_{i,k_i}\}$ for $i = 1, \ldots, n_v$.

For each family of models $\mathcal{G}_i$, we have a combined answer vector $\bold{p}_{\mathcal{G}_i} \in \R^{k_i}$ with the accuracy function $\Phi_{\mathcal{G}_i}$. Define $\Psi_{\mathcal{G}_i}$ \emph{the weight function for the family of models $G_i$} as
\[
\Psi_{\mathcal{G}_i} := 
\begin{cases}
	\frac{\Phi_{\mathcal{G}_i}}{\sum \Phi_{\mathcal{G}_i}} \quad \text{if } \sum \Phi_{\mathcal{G}_i} \neq 0, \\
	0 \quad \text{otherwise},
\end{cases}
\]
for each $i = 1,\ldots, n_v$.

Since each $\bold{p}_{\mathcal{G}_i}(\alpha,x)_{j}$, the $j$-th component of weighted answer from $\mathcal{G}_i$ for $x$ at certainty threshold $\alpha$ indicates how much it predicts $x$ to be classified into target $t_{i,j} \in \widetilde{T}_i$ at certainty threshold $\alpha$, we can multiply these `scores' to compute the overall agreement on $\cap_{i=1}^{n_v} t_{i,j_i}$ among the families $G_i$. For a given tuple $\left(t_{i,j_i}\right)_{i} \in \prod_{i=1}^{n_v}\widetilde{T}_{v,i}$ with multi-index $J = (j_1,\ldots, j_{n_v})$, define \emph{combined answer on $t_J$ (at node $v$)} as
\[
p_J(\alpha,x) = \prod_{i=1}^{n_v}\bold{p}_{\mathcal{G}_i}(\alpha,x)_{j_i},
\]
where $t_J$ be the tuple $\left(t_{i,j_i}\right)_{i=1,\ldots,n_v}$. Let $\mathcal{J}_v$ be the set of all possible indices $J$.

Let $\overline{T}_v = \{\bar{t}_1, \ldots, \bar{t}_k\}$ be the collection of refinements at node $v$. Then there exist unique indices $J_1,\ldots, J_k$ such that $\bar{t}_s = \bigcap t_{J_s} := t_{1,j_1}\cap \cdots \cap t_{n_v,j_{n_v}}$ where $J_s = (j_1,\ldots,j_{n_v})$ for each $s = 1,\ldots, k$. For each multi-index $J$ in $\mathcal{J}_v$, we call $J$ an \emph{invalid combination} if $\cap t_J = \emptyset$ and a \emph{valid combination} otherwise, that is, $\cap t_J \in \overline{T}_v$.
For an invalid combination $J = (j_1, \ldots, j_{n_v})$, we define \emph{the contribution factor} $\beta$ to distribute their combined answer $p_J$ to other refinements as follows:
\begin{equation}\label{eq:beta}
	\beta(t_{i,j_i},\bar{t}_s) = \begin{cases}
		\frac{p_{J_s}}{\sum_{\bar{t}_m \subset t_{i,j_i}}p_{J_m}} & \textrm{ if } \bar{t}_s \subset t_{i,j_i}\\
		0 & \textrm{ otherwise,}
	\end{cases}
\end{equation}
for each $i = 1,\ldots, n_v$ and a valid refinement $\bar{t}_s \in \overline{T}_v$.
Since $p_J = p_J \cdot \left(\Psi_{\mathcal{G}_1} + \cdots +\Psi_{\mathcal{G}_{n_v}} \right)$, we can decompose $p_J$ into 
\begin{align*}
	p_J = & ~ p_J \Psi_{\mathcal{G}_1} \beta(t_{1,j_1}, \bar{t}_1) + \cdots + p_J \Psi_{\mathcal{G}_{n_v}} \beta(t_{{n_v},j_{n_v}}, \bar{t}_1) \\
	& \quad + \cdots + p_J \Psi_{\mathcal{G}_1} \beta(t_{1,j_1}, \bar{t}_k) + \cdots + p_J \Psi_{\mathcal{G}_{n_v}} \beta(t_{{n_v},j_{n_v}}, \bar{t}_k),
\end{align*} as $\sum_{s = 1}^{k} \beta(t_{i,j_i},\bar{t}_s) = 1 $ for any $t_{i,j_i}$.
Then, we define $\mathcal{A}$ as \emph{the answer at node $v$}, a function from $X$ to $\R^k$ at certainty threshold $\alpha$, as
\[
\mathcal{A}_v = \left(p_{J_s} + \sum_{\substack{J \in \mathcal{J}_{v,\textrm{invalid}} \\ J = (j_1,\ldots,j_{n_v})}}\sum_{i=1}^{n_v}{p_J \Psi_{\mathcal{G}_i} \beta(t_{i,j_i}, \bar{t}_s) }\right)_{s = 1,\ldots, k} ,
\] where $\mathcal{J}_{v,\textrm{invalid}}$ is the collection of invalid combinations at node $v$. See Example~\ref{ex:refinement}

\begin{ex}[An example of voting procedure at a node]\label{ex:refinement}
	For a given node $v \in \mathfrak{T}$, let the flattened target $T_v = \{c_1,c_2,c_3,c_4, c_5\}$ and the indices of models $\mathcal{I}_v = \{(1,1),(1,2),2\}$ be associated with $v$. Suppose that the two models $g_{1,1}$ and $g_{1,2}$ share the target set $\widetilde{T}_1$, and $\widetilde{T}_2$ denotes the target set of $g_2$, where  
	\begin{align*}
		\widetilde{T}_1 &= \{t_{1,1}, t_{1,2}\} = \{\{c_1,c_2,c_3\}, \{c_4, c_5\}\}, \\
		\widetilde{T}_2 &= \{t_{2,1}, t_{2,2}, t_{2,3}\} = \{\{c_1,c_2\}, \{c_3,c_4\}, \{c_5\}\}.
	\end{align*}
	Their refinement $\overline{T}_v$ is $\{\{c_1,c_2\}, \{c_3\}, \{c_4\}, \{c_5\}\}$. Let $\mathcal{G}_1 = \{g_{1,1}, g_{1,2}\}$ and $\mathcal{G}_2 = \{g_2\}$, the collections of models sharing the same targets.
	
	Suppose that the accuracy functions $\Phi_{1,1}, \Phi_{1,2}, \Phi_2$ of models $g_{1,1}, g_{1,2}, g_2$ are given, respectively.
	For simplicity, we will look at certainty threshold $\alpha_0 = 0$ and suppress our notation reserved for the certainty threshold $\alpha$.
	Let an instance $x$ be given, and trained models provide answers for $x$ as follows:
		\begin{align*}
		g_{1,1} &= \left(a_{1,1}, a_{1,2}\right) & 
		g_{1,2} &= \left(a_{2,1}, a_{2,2}\right) & 
		g_{2} &= \left(b_1, b_2, b_3\right).
	\end{align*}
	
	Then we have the weighted answers $\bold{p}_\mathcal{G}$ from each collection of models sharing the same targets $\mathcal{G}$ :
	\begin{align*}
		\bold{p}_{\mathcal{G}_1} & =\left( \frac{\Phi_{1,1} a_{1,1} + \Phi_{1,2} a_{2,1}}{\Phi_{\mathcal{G}_1}} , \frac{\Phi_{1,1} a_{1,2} + \Phi_{1,2} a_{2,2}}{\Phi_{\mathcal{G}_1}} \right) := (a_1, a_2) ,\\
		\bold{p}_{\mathcal{G}_2} & =\left(b_1, b_2, b_3 \right),
	\end{align*}
	where $\Phi_{\mathcal{G}_1} = \Phi_{1,1}  + \Phi_{1,2}$ and $\Phi_{\mathcal{G}_1}  = \Phi_2$ are the accuracy functions of $\mathcal{G}_1$ and $\mathcal{G}_2$, respectively. Additionally, we have the weight functions $\Psi_{\mathcal{G}_1} = \frac{\Phi_{\mathcal{G}_1}}{\Phi_{\mathcal{G}_1} + \Phi_{\mathcal{G}_2}}$ and $\Psi_{\mathcal{G}_2} = \frac{\Phi_{\mathcal{G}_2}}{\Phi_{\mathcal{G}_1} + \Phi_{\mathcal{G}_2}}$.
	
	Let $\mathcal{J}_v$ the collection of all combinations be $\{J_1, J_2, J_3, J_4, J_5, J_6\}$ where
	\begin{align*}
		t_{J_1} &= \{t_{1,1}, t_{2,1}\}, ~ & t_{J_2} & = \{t_{1,1}, t_{2,2}\}, ~& t_{J_3} & = \{t_{1,1}, t_{2,3}\}, \\
		t_{J_4} & = \{t_{1,2}, t_{2,1}\}, ~&t_{J_5} & = \{t_{1,2}, t_{2,2}\}, ~&t_{J_6} & = \{t_{1,2}, t_{2,3}\}.
	\end{align*}
	Note that $J_3$ and $J_4$ are invalid combinations at node $v$, and the refinements of valid combinations are $\cap t_{J_1} = \{c_1,c_2\},~ \cap t_{J_2}  = \{c_3\},~ \cap t_{J_5}  = \{c_4\}$, and $ \cap t_{J_6}  = \{c_5\}$. We compute the combined answer for each $t_{J_i}$ as
	\begin{align*}
		p_{J_1} = & ~ a_1b_1,~&	p_{J_2} = & ~ a_1b_2,~&	p_{J_3} = &~  a_1b_3,~ \\
			p_{J_4} = &~  a_2b_1,~&	p_{J_5} = &~  a_2b_2,~&	p_{J_6} = &~  a_2b_3.
	\end{align*}
	Then the nontrivial contribution factors $\beta$ of $J_3$ defined in Equation~\ref{eq:beta} are
	\begin{align*}
		 \beta(t_{1,1}, \{c_1,c_2\}) &= \frac{p_{J_1}}{p_{J_1} + p_{J_2}} = \frac{b_1}{b_1+b_2} , ~ \\
		 \beta(t_{1,1}, \{c_3\}) &=  \frac{p_{J_2}}{p_{J_1} + p_{J_2}}  = \frac{b_2}{b_1+b_2}, ~ \\
		 \beta(t_{2,3}, \{c_5\}) &= 1, 
	\end{align*} as $t_{1,1} = \{c_1,c_2,c_3\}$ and $t_{2,3} = \{c_5\}$, and those $\beta$ of $J_4$ are
	\begin{align*}
		\beta(t_{1,2}, \{c_4\}) & = \frac{b_2}{b_2+b_3}, ~  \\
		\beta(t_{1,2}, \{c_5\}) & =  \frac{b_3}{b_2+b_3}, ~ \\
		\beta(t_{2,1}, \{c_1,c_2\}) & = 1, 
	\end{align*} as $t_{1,2} = \{c_4,c_5\}$ and $t_{2,1} = \{c_1,c_2\}$.
    Therefore, the answer $\mathcal{A}_v$ at node $v$ for $x$ is
    \begin{align*}
    \mathcal{A}_v(x) & =  \left(a_1b_1 + \frac{a_1b_1b_3 \Psi_{\mathcal{G}_1}}{b_1 + b_2} + a_2b_1 \Psi_{\mathcal{G}_2} ,~ a_1b_2 + \frac{a_1b_2b_3\Psi_{\mathcal{G}_1}}{b_2+b_3} \right.\\
    & \quad , \left. a_2b_2 + \frac{a_2b_1 b_2\Psi_{\mathcal{G}_1}}{b_2 + b_3}, ~a_2b_3 + a_1 b_3\Psi_{\mathcal{G}_2} + \frac{a_2 b_1b_3\Psi_{\mathcal{G}_1}}{b_2+b_3} \right).
    \end{align*} \end{ex}

\subsubsection{Accumulation of votes along valid paths}



For a target class $c \in T$, a sequence of nodes $\gamma = (s_0, s_1, \ldots, s_m)$ in $\mathfrak{T}$ is called  \emph{a valid path for $c$} if $\gamma$ satisfies the following conditions:
\begin{enumerate}
	\item $s_0$ is the root of $\mathfrak{T}$.
	\item $c \in T_{s_i}$ for all $i = 0,\ldots,m$.
	\item $\overline{T}_{s_m}$ consists of thin targets.
\end{enumerate}
Let $\gamma = (s_0, s_1, \ldots, s_m)$ be a valid path for a class $c \in T$, where $s_0$ is the root of $\mathfrak{T}$.
Since each trained model provides prediction independently, we define $\cM(\gamma,\alpha,x)$ \emph{the weighted answer for $x$ at certainty threshold $\alpha$ along a valid path $\gamma$} as the product \begin{equation}\label{eqn:accumulationVotesAlongPath}
	\cM(\gamma,\alpha,x) = \left(\prod_{i=0}^{m}A_{s_i}(\alpha,x)_{j_{i,t}}\right)_{t \in \overline{T}_{s_m}} = \left(\cM_{t_1},\ldots,\cM_{t_{|\overline{T}_{s_m}|}}\right),
\end{equation} which represents how much $\cM$ predicts $x$ to be classified in each target $t \in \overline{T}_{s_m}$ along the path $\gamma$. Here, $j_{i,t}$ is the index of $\bar{t} \in \overline{T}_{s_i}$ such that $\bar{t}$ is the the unique refinement in $\overline{T}_{s_i}$ containing $c$. Then define $\bold{P}(\gamma,\alpha,x)$ \emph{the prediction for $x$ at certainty threshold $\alpha$ along a valid path $\gamma$} as the $\arg\max \cM(\gamma,\alpha,x)$.
\begin{remark}
Under the specialization method explained in the Section \ref{specialization}, we can construct `gating layer' as in the Mixture of Expert\cite{AdaptiveMixtureofLocalExperts} using this voting strategy.
Let $g$ be a trained model in $\cM$ and the targets of $g$ be $\widetilde{T} = \left\{ t_1,\ldots,t_k\right\}$ where $t_i = \{c_{i,1},\ldots,c_{i,n_i}\}$ for $i=1,\ldots,k$ such that $T= \{c_{1,1},\ldots,c_{k,n_k}\}$. Then $g$ serves as the `gating' layer in $\cM$ navigating an instance to other trained models that are trained on dataset containing classes exclusively within a target $t_i$ of $\widetilde{T}$. See \cite{Jung-Lau} for the experimental results.
\end{remark}
\subsubsection{Vote using validation history}
We introduce the `using validation history' method in prediction to alleviate concerns regarding the optimal valid path or certainty threshold. In other words, $\alpha$ and $\gamma$ in the Equation \ref{eqn:accumulationVotesAlongPath} are fixed through this method based on the validation dataset.


Let $X_{\textrm{val}}$ be a measurable subset of $X$ with $\mu(X_{\textrm{val}})<\infty$, where $\mu$ is the measure of $X$. Let $\Gamma(c)$ denote the set of all valid paths in the target graph $\mathfrak{T}$ for a class $c$. Given the true label function $f : X \to T$, we define $r$ as \emph{the expected accuracy along a path $\gamma$ at certainty threshold $\alpha$ to class $c$} as:
\[
r_c(\gamma,\alpha) := \frac{\left|\left(\{\bold{P}(\gamma, \alpha, x) = c\}\cap f^{-1}(c)\right) \cup \{\bold{P}(\gamma, \alpha, x) \neq c \textrm{ and } f(x) \neq c\} \right|}{\left|X_\textrm{val} \right|},
\] where $\gamma \in \Gamma(c)$.
Since there are finite number of valid paths for each target class and $\alpha \in [0,1]$, there exists a tuple of maximizers $(\gamma^*(c), \alpha^*(c))$ for each class $c$ such that the expected accuracy $r$ attains its supremum at $\alpha^*(c)$. We define \emph{the answer using validation history for $x \in X$} as $$\cM(x) = \left(\cM_{c_1}(\gamma^*(c_1), \alpha^*(c_1),x),\ldots,\cM_{c_N}(\gamma^*(c_N), \alpha^*(c_N),x) \right).$$

\begin{remark}
	Let a system of trained models $\cM$ and validation dataset $X_{\textrm{val}}$ be given.
	The combining rule using validation history for each trained model serves as the role of $\rho$ in the Definition \ref{def:fuzlogifold}, and $P : X \to [0,1]$ is defined by the vote of $\cM$ using the validation history.
\end{remark}
%

\newpage
\appendix
\section{Pseudo-Algorithms}\label{appendixA}
In this appendix, we provide pseudocode implementations for the algorithms discussed in Section~\ref{sec:algo}. Given classification classes are $c_1, \ldots ,c_N$.

\begin{algorithm}
	\caption{Refinement}\label{alg:Refinement}
	\begin{algorithmic}[1]
		\Require Targets $T_1, \ldots, T_n$ with $T = \cup T_1 = \cdots = \cup T_n$ 
		\Ensure Refinement $\overline{T}$
		\State Initialize $\overline{T}$
		\State Combinations $\gets$ all combinations  $(t_1,\ldots,t_n) \in \prod_i^n T_i$
		\For{$c = (t_1,\ldots,t_n) \in$ Combinations}
		\State $\bar{t} \gets \bigcap c = t_1 \cap \cdots \cap t_n$
		\State $c.\text{refinement} \gets \bar{t}$
		\If{$\bar{t} \neq \emptyset$}
		\State $\bar{t}.\text{component} \gets c$
		\State Append $\bar{t}$ to $\overline{T}.\text{valid}$
		\State Append $c$ to $\overline{T}.\text{validCombinations}$
		\EndIf
		\State Append $c$ to $\overline{T}.\text{allCombinations}$
		\EndFor
		\State \Return $\overline{T}$
	\end{algorithmic}
\end{algorithm}

For instance, in Algorithm~\ref{alg:Refinement} with Example~\ref{ex:refinement}, we have six combinations: \begin{align*}
	\textrm{CB}_1 :=& (t_{1,1}, t_{2,1}) = \left(\{c_1,c_2,c_3\},\{c_1,c_2\}\right), & \textrm{CB}_4 := (t_{1,2}, t_{2,1})& = \left(\{c_4,c_5\},\{c_1,c_2\}\right),\\ 
	\textrm{CB}_2 := &(t_{1,1}, t_{2,2}) = \left(\{c_1,c_2,c_3\},\{c_3,c_4\}\right), & \textrm{CB}_5 := (t_{1,2}, t_{2,2})& = \left(\{c_4,c_5\},\{c_3,c_4\}\right),\\
	\textrm{CB}_3 := &(t_{1,1}, t_{2,3}) = \left(\{c_1,c_2,c_3\},\{c_5\}\right),& \textrm{CB}_6 := (t_{1,2}, t_{2,3})& = \left(\{c_4,c_5\},\{c_5\}\right). 	
\end{align*}

Let $\bar{t}_i$ denote the refinement obtained by the combination $\textrm{CB}_i$ for $i=1,\ldots,6$. In other words, $\textrm{CB}_i.\textrm{refinement} = \bar{t}_i$ and $\bar{t}_i.\textrm{component} = \textrm{CB}_i$. Then $\bar{t}_1 = \{c_1,c_2\}, \bar{t}_2 = \{c_3\}, \bar{t}_3 = \bar{t}_4 = \emptyset, \bar{t}_5 = \{c_4\}$, and $\bar{t}_6 = \{c_5\}$. $\textrm{CB}_3$ and $\textrm{CB}_4$ are invalid combinations, and $\bar{t}_1,\bar{t}_2, \bar{t}_5$ and $\bar{t}_6$ are valid refinements. Therefore, we have 
\begin{align*}
	\overline{T}.\textrm{allCombinations} =& \left\{\textrm{CB}_i : i= 1,\ldots, 6\right\},\\ \overline{T}.\textrm{validCombinations} =& \left\{\textrm{CB}_1,\textrm{CB}_2,\textrm{CB}_5,\textrm{CB}_6\right\}, \\ \overline{T}.\textrm{valid}=&\{\bar{t}_1,\bar{t}_2,\bar{t}_5,\bar{t}_6\}.
\end{align*}

\begin{algorithm}
	\caption{Construct Target Graph}\label{alg:Target Graph}
	\begin{algorithmic}[1]
		\Require Trained models $g_1, \ldots, g_n$ and corresponding targets $\widetilde{T}_1,\ldots,\widetilde{T}_n$.
		\Ensure Target graph $\mathfrak{T}$.
		\State Re-index $\{\widetilde{T}_1,\ldots,\widetilde{T}_n\}$ as $\{\widetilde{T}_1,\ldots,\widetilde{T}_m\}$ such that all elements are distinct.
		\State $G_1, \ldots ,G_m \gets $ corresponding collections of trained models associated with $\widetilde{T}_1, \ldots, \widetilde{T}_m$ \Comment{Group models sharing the same targets together}
		\State Initialize an array $\mathcal{T}$
		\For{$\widetilde{T}$ runs over $ \widetilde{T}_1,\ldots,\widetilde{T}_m$} 
		\State $T \gets \bigcup\widetilde{T}$
		\If{$T \notin \mathcal{T}$}  Append $T$ to $\mathcal{T}$
		\EndIf
		\EndFor
		
		\State Sort $\mathcal{T}$ by decreasing size
		\For{$T$ runs over $\mathcal{T}$}
		\State Add node $T$ to $\mathfrak{T}$
		\State Initialize $T.\textrm{nextNodes}$
		\State $i \gets$ be the index of $T$ in $\mathcal{T}$
		\For{$T' \in \mathcal{T}[0,\ldots,i-1]$}
		\If{$T \subset T'$}  Append $T$ to  $T'.\textrm{nextNodes}$ \EndIf
		\EndFor
		
		\State Initialize $T.\textrm{targetsAndModels}$
		\For{$\widetilde{T} \in \{\widetilde{T}_1,\ldots,\widetilde{T}_m\}$}
		\If{$\cup \widetilde{T} = T$}
		\State Add $\widetilde{T}$ and corresponding group of models $G$ to $T.\textrm{targetsAndModels}$
		\EndIf
		\EndFor
		\State Initialize $T.\textrm{models}$
		\State $\left\{\left(\widetilde{T}_j, G_j\right)_{j \in \Lambda}\right\} \gets T.\textrm{targetsAndModels}$ \Comment{$\Lambda$ is a finite index set}
		
		\State Add $\bigcup_{j\in \Lambda} G_j $ to $T.\textrm{models}$
		\State $T.\textrm{refinement} \gets \Call{Refinement}{\left(\widetilde{T}_j\right)_{j \in \Lambda}}$
		
		\EndFor

		\State\Return $\mathfrak{T}$
		
	\end{algorithmic}
\end{algorithm}

\begin{algorithm}
	\caption{Fuzzy Accuracy}\label{alg:FuzzyAccuracy}
	\begin{algorithmic}[1]
		\Require Threshold $\alpha\in[0,1]$, Model $g$, Dataset $\mathcal{Z} =\{(x,y_x)\}$
		\Ensure Accuracy $\Phi$ and Certain Part
		\State Initialize $\text{CertainPart}$ and $\Phi$
		\State $\text{CertainPart}(\alpha) \gets \{x \mid \max g(x) \geq \alpha\}$
		\State $\Phi(\alpha) \gets \frac{|\{x \in \text{CertainPart}(\alpha) \mid \arg\max g(x) = y_x\}|}{|\text{CertainPart}(\alpha)|}$
		\State \Return $\Phi,\text{CertainPart}$
	\end{algorithmic}
\end{algorithm}
\begin{algorithm}
	\caption{Total Weight and Rho}\label{alg:WeightByAccuracy}
	\begin{algorithmic}[1]
		\Require Threshold $\alpha\in[0,1]$, Models $G= \{g_1, \ldots, g_n\}$ with accuracies $\Phi_{g_1}, \ldots, \Phi_{g_n}$,  Instance $x$
		\Ensure Sum of weight by Accuracy $\Phi$ and $\left(\rho_g\right)_{g\in G}$
		\For{$g \in G$}
		\If{$\max g(x) \geq \alpha$}
		\State $\rho_g (\alpha,x)\gets 1$
		\Else
		\State $\rho_g(\alpha,x)\gets \epsilon$ \Comment{$\epsilon$ can be any sufficiently small number.}
		\EndIf
		\EndFor
		\State $\Phi(\alpha,x) \gets \sum_{g \in G} \rho_g(\alpha,x)\Phi_g(\alpha)$
		\State\Return $\Phi, \left(\rho_g\right)_{g\in G}$
	\end{algorithmic}
\end{algorithm}

\begin{algorithm}
	\caption{Voting Rule for Shared Targets}\label{alg:VotingSharedTargets}
	\begin{algorithmic}[1]
		\Require Threshold $\alpha\in[0,1]$, Targets $T=\{t_1, \ldots , t_n\}$, Models $G= \{g_1, \ldots, g_k\}$ with accuracies $\Phi_{g_1}, \ldots, \Phi_{g_k}$, Instance $x$
		\Ensure Answer $\bold{p}$ 
		\State $\Phi_G, \left(\rho_g\right)_{g\in G} \gets \Call{Total Weight and Rho}{\alpha, G, x}$
		\State $\bold{p}(\alpha,x) \gets \left(\Phi^{-1}(\alpha,x) \cdot \sum_{g \in G} \rho_g(\alpha,x) \Phi_g(\alpha,x) g_j(x)\right)_{j=1,\ldots,n}$
		\State \Return $\bold{p}$
	\end{algorithmic}
\end{algorithm}

\begin{algorithm}
	\caption{Distribute Answers at a node}\label{alg:Distribute Answers}
	\begin{algorithmic}[1]
		
		\Require Node $v$ in a target graph $\mathfrak{T}$ with $\overline{T} = v.\textrm{refinement}$, Combined votes $\{p_c\}_{c \in \overline{T}.\textrm{allCombinations}}$, Instance $x$
		\Ensure Answer for $x$ at each threshold $\alpha$
		\State $\left\{\left(\widetilde{T}_i,G_i\right)_{i=1,\ldots,k}\right\} \gets v.\textrm{targetsAndModels}$
		\State Initialize $\mathfrak{p}$ and $\beta$
		\For{$t \in \bigcup_{i=1}^{k}\widetilde{T}_i$}
		\State $\overline{T}_{t} \gets \{ \bar{t} \in \overline{T}.\textrm{valid} : \bar{t} \subset t \}$
		\State $C_t \gets \{\bar{t}.\textrm{component} : \bar{t}\in \overline{T}_t\}$
		\EndFor
		\For{$t \in \bigcup_{i=1}^{k}\widetilde{T}_i$}
		\For{$\bar{t} \in \overline{T}_{t} $}
		\State $c \gets \bar{t}.\textrm{component}$
		
		\State $\beta(t,\bar{t}) \gets \frac{p_c(\alpha,x)}{\sum_{d \in C_t} p_d(\alpha,x)}$
		\EndFor
		\EndFor
		
		\For{$c=\left(t_{1,j_1},\ldots,t_{k,j_k}\right) \in \overline{T}.\textrm{allCombinations}$}
		\If{$c \in \overline{T}.\textrm{validCombinations}$}
		\State $s \gets $ the index of $c.\textrm{refinement}$ in $\overline{T}$.
		\State $\mathfrak{p}(\alpha,x,c)_s \gets p_c(\alpha,x)$
		\State $\mathfrak{p}(\alpha,x,c)_j \gets 0$ for all $j \neq s$
		\Else
		
		\For{$s=1,\ldots,m$}
		\If{$\{t \in c : \bar{t}_s \subset t\}\neq \emptyset$}
		\State $\mathfrak{p}(\alpha,x,c)_s \gets \sum_{t_{i,j_i} \in \{t \in c : \bar{t}_s \subset t\}} \beta(t_{i,j_i}, \bar{t}_s) \Psi_{G_i}(\alpha,x) p_c(\alpha,x)$
		\Else
		\State $\mathfrak{p}(\alpha,x,c)_s \gets \epsilon $ \Comment{$\epsilon$ can be any sufficiently small number.}
		\EndIf
		\EndFor
		
		\EndIf
		\EndFor
		\State\Return $\mathfrak{p}$
	\end{algorithmic}
\end{algorithm}

\begin{algorithm}
	\caption{Voting rule at a node}\label{alg:Voting rule at a node}
	\begin{algorithmic}[1]
		
		\Require Threshold $\alpha\in[0,1]$, Node $v$ in a target graph $\mathfrak{T}$, Instance $x$
		\Ensure Answer for $x$ at each threshold $\alpha$

		\State $\left\{\left(\widetilde{T}_i,G_i\right)_{i=1,\ldots,k}\right\} \gets v.\textrm{targetsAndModels}$
		\State $\mathcal{G} \gets \{G_1,\ldots,G_k\}$
		
			\For{$(\widetilde{T}, G) \in v.\textrm{targetsAndModels}$}
			\State $\bold{p}_{G} \gets \Call{Voting Rule for Shared Targets}{\alpha,\widetilde{T},G}$
			\State $\Phi, \left(\rho_g\right)_{g\in G} \gets \Call{Total Weight and Rho}{\alpha, G, x}$
			\EndFor
			
		\State $\Psi \gets \left(\frac{\Phi_G}{\sum_{G \in \mathcal{G}}\Phi_G}\right)_{G \in \mathcal{G}} $

		\State $\overline{T} = \{\bar{t}_1,\ldots,\bar{t}_m\}  \gets v.\textrm{refinement}$

			\For{$t \in \bigcup_{i=1}^{k}\widetilde{T}_i$}
			\State $\overline{T}_{t} \gets \{ \bar{t} \in \overline{T}.\textrm{valid} : \bar{t} \subset t \}$
			\State $C_t \gets \{\bar{t}.\textrm{component} : \bar{t}\in \overline{T}_t\}$
			\EndFor

		\State Enumerate $\{\widetilde{T}_1,\ldots,\widetilde{T}_k\} = \{\{t_{1,1},\ldots, t_{1,r_1}\},\ldots, \{t_{k,1},\ldots, t_{k,r_k}\}\}$
		\State $C \gets \overline{T}.\textrm{allCombinations}$
			
				\For{$c=\left(t_{1,j_1},\ldots,t_{k,j_k}\right) \in  C$}
				    \State Initialize $p_c$
				    \State $p_c(\alpha,x) \gets \prod_{i=1}^{k}\bold{p}_{G_i}(\alpha,x)_{j_i}$
    			\EndFor
    
   		\State $\textrm{Answers} \gets \Call{Distribute Answers at a node}{v,\{p_c\}_{c \in C}, x}$
		\State $\mathfrak{p}(\alpha,x) \gets \left(\sum_{c \in C} \textrm{Answers}(\alpha,x,c)_j\right)_{j=1,\ldots,m}$
		   	
   		\State\Return $\mathfrak{p}$
	\end{algorithmic}
\end{algorithm}

\begin{algorithm}
	\caption{Valid Paths}\label{alg:Valid Paths}
	\begin{algorithmic}[1]
		
		\Require Target Graph $\mathfrak{T}$, Target class $c \in \{c_1,\ldots, c_N\}$
		\Ensure $\textrm{Valid Paths}$
		\State $s_0 \gets \{c_1,\ldots, c_N\}$
		\State Initialize $\textrm{Candidates}$
		\State Search for all paths $s_0, \ldots, s_m$ such that $s_i \in s_{i-1}.\textrm{nextNodes}$ and append to $\textrm{Candidates}$
		\State Initialize $\textrm{Valid Paths}$
			\For{$s_0,\ldots,s_m \in \textrm{Candidates}$}
			\If{$c \in s_i$ for all $i=0,\ldots,m$ and $|\bar{t}| = 1 \textrm{ for all } \bar{t}\in s_m.\textrm{refinement}$}
			\State Append $s_0,\ldots,s_m$ to $\textrm{Valid Paths}$
			\EndIf
			\EndFor
		\State\Return $\textrm{Valid Paths}$
		
	\end{algorithmic}
\end{algorithm}

\begin{algorithm}
	\caption{Voting rule along a path}\label{alg:Voting rule along a path}
	\begin{algorithmic}[1]
		
		\Require Threshold $\alpha\in[0,1]$, Valid path $\gamma = \left(s_0,\ldots,s_m\right)$ to $c$, Instance $x$
		\Ensure $\cM$ Answer for $x$ at each threshold $\alpha$ along a valid path $\gamma$ 
		\State $\{\bar{t}_1,\ldots,\bar{t}_n\} \gets s_m.\textrm{refinement}$
		
			\For{$\bar{t}_j \in \{\bar{t}_1,\ldots,\bar{t}_n\} $}
				\For{$s_i \in \gamma$}
				\State $\textrm{Answer}_{s_i}(\alpha,x)\gets \Call{Voting rule at a node}{\alpha,s_i,x}_k$ where $\bar{t}_k \in s_i.\textrm{refinement}$ contains $c$
				\EndFor
			\State $\cM(\gamma, \alpha, x) \gets \prod_{i=0}^{m} \textrm{Answer}_{s_i}(\alpha,x)$
			\EndFor
		\State\Return $\cM$
	\end{algorithmic}
\end{algorithm}

\begin{algorithm}
	\caption{Compute Expected Accuracy}\label{alg:Compute Expected Accuracy}
	\begin{algorithmic}[1]
		
		\Require Answer Vector $\cM$, Target Class $c$, Valid Path $\gamma$ to $c$, Threshold $\alpha \in [0,1]$, (Validation) Dataset $X = \{(x,y_x)\}$
		\Ensure Expected Accuracy $r$
		\State $\textrm{Prediction}(\gamma,\alpha,x) \gets \arg\max \cM(\gamma,\alpha,x)$
		\State $X_c \gets \{x \in X: y_x = c\}$
		\State $TP \gets \left\{\textrm{Prediction}(\gamma,\alpha,x) = c\right\} \cap X_c$
		\State $TN \gets \left\{\textrm{Prediction}(\gamma,\alpha,x) \neq c\right\} \cap X\setminus X_c$
		\State $r(\gamma,\alpha) \gets \frac{\left|TP \cup TN\right|}{\left|X\right|}$
		\State \Return $r$
	\end{algorithmic}
\end{algorithm}

\begin{algorithm}
	\caption{Vote using Validation history}\label{alg:Vote using Validation history}
	\begin{algorithmic}[1]
		
		\Require Thresholds $A=\{0 = \alpha_0, \ldots, \alpha_n\}$, Target Graph $\mathfrak{T}$, (Validation) Dataset $\mathcal{Z} = \{(x,y_x)_{x \in X_{\textrm{val}}}\}$, Instance $x$
		\Ensure Final Answer $\cM(x) = \left(\cM_{c_1}(x),\ldots, \cM_{c_N}(x)\right)$
			\For{$c \in \{c_1,\ldots,c_N\}$}
			\State $\Gamma(c) \gets \Call{Valid Paths}{\mathfrak{T}, c}$
				\For{$\gamma=(s_0,\ldots,s_m) \in \Gamma(c)$, and $\alpha \in A$}
					\State $\overline{T} \gets s_m.\textrm{refinement}$
					\For{$x \in \mathcal{X}_{\textrm{val}}$ and $\bar{t} \in \overline{T}$}
					\State $\cM_{\bar{t}}(\gamma,\alpha,x) \gets \Call{Voting rule along a path}{\alpha,t, \gamma,x}$
					\EndFor
					\State $\cM(\gamma,\alpha,x) \gets \left(\cM_{\bar{t}}\right)_{\bar{t} \in \overline{T}}$
					\State $r(c,\gamma,\alpha) \gets \Call{Compute Expected Accuracy}{\cM, c,\gamma,\alpha,X_\textrm{val}}$
				\EndFor
				\State $\gamma^*(c), \alpha^*(c) \gets \arg\max_{\gamma, \alpha}r(c,\gamma,\alpha)$
				\State $\cM_c(x) \gets \cM(\gamma^*(c),\alpha^*(c),x)$
			\EndFor
		\State \Return $\cM_{c_1}(x),\ldots,\cM_{c_N}(x)$
	\end{algorithmic}
\end{algorithm}
%


\clearpage
	\bibliographystyle{siam}
	\bibliography{geometry}{}
\end{document}